\documentclass[12pt]{article}

\usepackage{latexsym,amsfonts,amsmath,amssymb,graphicx,hyperref}

\newcommand{\comment}[1]{}

\newcommand{\Om}{\Omega}

\newcommand{\EE}{\mbox{\bf E}\,}

\newcommand{\R}{\mathbb{R}}
\newcommand{\C}{\mathbb{C}}
\newcommand{\Q}{\mathbb{Q}}
\newcommand{\HH}{\mathbb{H}}
\newcommand{\N}{\mathbb{N}}

\newcommand{\Z}{\mathbb{Z}}

\newcommand{\pa}{\partial}

\newcommand{\F}{{\cal F}}

\newcommand{\no}{\noindent}

\newcommand{\BGE}{\begin{equation}}
\newcommand{\BGEN}{\begin{equation*}}
\newcommand{\EDE}{\end{equation}}
\newcommand{\EDEN}{\end{equation*}}

\def\eps{\varepsilon}
\def\til{\widetilde}
\def\ha{\widehat}
\def\sem{\setminus}
\def\lin{\overline}
\def\vphi{\varphi}

\def\del{\delta}

\DeclareMathOperator{\hcap}{hcap} \DeclareMathOperator{\id}{id}
\DeclareMathOperator{\Imm}{Im } 
 
\DeclareMathOperator{\SLE}{SLE} \DeclareMathOperator{\HP}{HP}

\def\h0{{\bf h}}
\def\luto{\stackrel{\rm l.u.}{\longrightarrow}}

\newtheorem{Lemma}{Lemma}[section]
\newtheorem{Theorem}{Theorem}[section]

\newtheorem{Proposition}{Proposition}[section]
\numberwithin{equation}{section}

\begin{document}
\title{\bf Reversibility of Chordal SLE}
\date{\today}
\author{Dapeng Zhan}
\maketitle
\begin{abstract}
We prove that the chordal SLE$_\kappa$ trace is reversible for
$\kappa\in(0,4]$.
\end{abstract}

\section{Introduction}
Stochastic Loewner evolutions (SLEs) are introduced by Oded Schramm
(\cite{S-SLE}) to describe the scaling limits of some lattice
models, whose scaling limits satisfy conformal invariance and Markov
property. The basic properties of SLE are studied in
\cite{RS-basic}. There are several different versions of SLE. A
chordal SLE is defined in a simply connected domain, which is about
some random curve in the domain that grows from one boundary point
to another.

So far it has been proved that chordal SLE$_6$ is the scaling limit
of explorer line of the site percolation on triangular lattice with
half open and half closed boundary conditions (\cite{SS-6} and
\cite{Conv-6}); chordal SLE$_8$ is the scaling limit of UST Peano
curve with half free and half wired boundary conditions
(\cite{LSW-2}); chordal SLE$_4$ is the scaling limit of contour line
of the two-dimensional discrete Gaussian free field with suitable
boundary values (\cite{SS}); and chordal SLE$_2$ is the scaling
limit of LERW started near one boundary point, conditioned to leave
the domain near the other boundary point (\cite{LERW}). In
\cite{LSW-8/3}, SLE$_{8/3}$ is proved to satisfy the restriction
property. From these results, we know that the chordal SLE$_\kappa$
trace is reversible for $\kappa=6,8,4,2,8/3$.

In \cite{RS-basic}, it is conjectured that the chordal SLE$_\kappa$
trace is reversible for all $\kappa\in[0,8]$. Oded Schramm and
Steffen Rohde are planning to prove that the chordal SLE$_\kappa$
trace is not reversible for $\kappa>8$ (\cite{Plan}). Scott
Sheffield proposed that the reversibility can be derived from the
relationship with the Gaussian free field (\cite{ICM}). In this
paper we will prove this conjecture for $\kappa\in(0,4]$ using only
techniques of probability theory and stochastic processes. The main
idea of this paper is as follows.

Suppose $(\beta(t))$ is a chordal SLE$_\kappa$ trace in a simply
connected domain $D$ from a prime end $a$ to another prime end $b$.
From Markov property of SLE, for a fixed time $t_0$, conditioned on
the curve $\beta([0,t_0])$, the rest of the curve $(\beta(t):t\ge
t_0)$ has the same distribution as a chordal SLE$_\kappa$ trace in
$D_{t_0}:=D\sem\beta([0,t_0])$ from $\beta(t_0)$ to $b$. Assume that
the chordal SLE$_\kappa$ trace is reversible. Then the reversal of
$(\beta(t):t\ge t_0)$ has the same distribution as the chordal
SLE$_\kappa$ trace in $D_{t_0}$ from $b$ to $\beta(t_0)$. On the
other hand, since $(\beta(t):t\ge t_0)$ is a part of the
SLE$_\kappa$ trace in $D$ from $a$ to $b$, so from the
reversibility, the reversal of $(\beta(t):t\ge t_0)$ should be a
part of SLE$_\kappa$ trace in $D$ from $b$ to $a$. Suppose $\gamma$
is an SLE$_\kappa$ trace in $D_{t_0}$ from $b$ to $\beta(t_0)$. From
the above discussion, if we integrate $\gamma$ against all possible
curves $\beta([0,t_0])$, we should get a part of SLE$_\kappa$ trace
in $D$ from $b$ to $a$, assuming that the chordal SLE$_\kappa$ trace
is reversible.

To prove the reversibility, we want to find a coupling of two
SLE$_\kappa$ traces, one is from $a$ to $b$, the other is from $b$
to $a$, such that the two curves visit the same set of points. If
such coupling exists, we choose a pair of disjoint hulls, each of
which contains some neighborhood of $a$ or $b$ in $\HH$, and stop
the two traces when they leave one of the two hulls, respectively.
Before these stopping times, the two traces are disjoint from each
other. The joint distribution of the two traces up to these stopping
time should agree with that of $\beta$ and $\gamma$ discussed in the
last paragraph up to the same stopping times. The Girsanov Theorem
suggests that this distribution is absolutely continuous w.r.t.\
that of two independent chordal SLE$_\kappa$ traces (one from $a$ to
$b$, the other from $b$ to $a$) stopped on leaving the above two
hulls. And the Radon-Nikodym derivative is described by a
two-dimensional local martingale, which has the property that when
one variable is fixed, it is a local martingale in the other
variable. This is the $M(\cdot,\cdot)$ in Theorem \ref{M}. It is
closely related with Julien Dub\'edat's work about commutation
relations for SLE (\cite{Julien}).

Using the $M(\cdot,\cdot)$, we may construct a portion of the
coupling up to certain stopping times. To construct the global
coupling, the difficulty arises when the two hulls collide, and the
absolute continuity blows up after that time. In fact, we can not
expect that the global coupling we are looking for is absolutely
continuous w.r.t.\ two independent SLE. Instead, the coupling
measure will be the weak limit of a sequence of absolutely
continuous coupling measures. Each measure in the sequence is
generated from some two-dimensional bounded martingale, which is the
$M_*(\cdot,\cdot)$ in Theorem \ref{martg}. The important property of
$M_*$ is that on the one hand, it carries the information of $M$ as
much as we want; on the other hand, it is uniformly bounded, and
remains to be a martingale even after the two hulls collide. So
$M_*$ can be used as the Radon-Nikodym derivative to define a global
coupling measure.

Although the results in this paper about martingales hold for all
$\kappa>0$, there are still some work to do to prove the
reversibility when $\kappa\in (4,8]$. The case that $\kappa=6$
illustrates the difficulty of this problem. In this case all
martingales are constant $1$. But the reversibility of SLE$_6$ is
still a non-trivial problem.

As suggested by Oded Schramm, the method in this paper can be used
to prove the duality of SLE. The author is planning to write another
paper about this object. We may use this to study the reversal of
the chordal SLE$(\kappa,\rho)$ trace (\cite{LSW-8/3}), continuous
LERW curve (\cite{LERW}), and annulus SLE trace
(\cite{Zhan}\cite{ann-prop}).

This paper is organized in the following way. In Section
\ref{prelim} we give the definition of chordal SLE and some other
basic notations, and then present the main theorem of this paper. In
Section \ref{ensem}, we study the relations of two SLE that grow in
the same domain. In Section \ref{Two-Dim} we present the
two-dimensional local martingale $M$, and check its property by
direct calculation of stochastic analysis. In Section
\ref{sec-bound}, we give some stopping times up to which $M$ is
bounded. And at the end of Section \ref{sec-bound}, we give a
detailed explanation of the meaning of $M$. In Section \ref{Const},
we use the local martingale to construct some two-dimensional
bounded martingale $M_*$. In Section \ref{coup}, we use $M_*$ to
construct a sequence of coupling measures. The limit of these
measures in some suitable sense is also a coupling measure. We
finally prove that under the limit measure, the two SLE$_\kappa$
traces coincide with each other.

\vskip 3mm

\no{\bf Acknowledgements}. I would like to thank Nikolai Makarov for
introducing me to the area of SLE. I also thank Oded Schramm for
some important suggestions about this paper and future work.

\section{Chordal Loewner Equation and Chordal SLE}\label{prelim}
Let $\HH=\{z\in\C:\Imm z>0\}$ denote the upper half complex plane.
If $H$ is a bounded closed subset of $\HH$ such that $\HH\sem H$ is
simply connected, then we call $H$ a hull in $\HH$ w.r.t.\ $\infty$.
For such $H$ there is a unique $\vphi_H$ that maps $\HH\sem H$
conformally onto $\HH$ such that $\vphi_H(z)=z+\frac{c}{z}+O(1/z^2)$
as $z\to\infty$ for some $c\ge 0$. Such $c$ is called the half-plane
capacity of $H$, and is denoted by $\hcap(H)$.

\begin{Proposition} Suppose $\Om$ is an open neighborhood of $x_0\in\R$
in $\HH$. Suppose $W$ maps $\Om$ conformally into $\HH$ such that
for some $r>0$, if $z\to(x_0-r,x_0+r)$ in $\Om$ then $W(z)\to\R$. So
$W$ extends conformally across $(x_0-r,x_0+r)$ by Schwarz reflection
principle. Then for any $\eps>0$, there is some $\del>0$ such that
if a hull $H$ in $\HH$ w.r.t.\ $\infty$ is contained in
$\{z\in\HH:|z-x_0|<\del\}$, then $W(H)$ is also a hull in $\HH$
w.r.t.\ $\infty$, and
$$|\hcap(W(H))-W'(x_0)^2\hcap(H)|\le\eps|\hcap(H)|.$$
\label{hcap}
\end{Proposition}\vskip -9mm
{\bf Proof.} This is Lemma 2.8 in \cite{LSW1}. $\Box$

\vskip 4mm

For a real interval $I$, let $C(I)$ denote the real valued
continuous function on $I$. Suppose $\xi\in C([0,T))$ for some
$T\in(0,+\infty]$. The chordal Loewner equation driven by $\xi$ is
as follows: \BGE
\pa_t\vphi(t,z)=\frac{2}{\vphi(t,z)-\xi(t)},\quad\vphi(0,z)=z.\label{chordal-equation}\EDE
For $0\le t<T$, let $K(t)$ be the set of $z\in\HH$ such that the
solution $\vphi(s,z)$ blows up before or at time $t$. We call $K(t)$
and $\vphi(t,\cdot)$, $0\le t<T$, chordal Loewner hulls and maps,
respectively, driven by $\xi$. Then for each $t\in [0,T)$,
$\vphi(t,\cdot)$ maps $\HH\sem K(t)$ conformally onto $\HH$. Suppose
for every $t\in[0,T)$,
$$\beta(t):=\lim_{z\in\HH,z\to \xi(t)}\vphi(t,\cdot)^{-1}(z)\in\HH\cup\R$$
exists, and $\beta(t)$, $0\le t<T$, is a continuous curve. Then for
every $t\in [0,T)$, $K(t)$ is the complement of the unbounded
component of $\HH\sem\beta((0,t])$. We call $\beta$ the chordal
Loewner trace driven by $\xi$. In general, such trace may not exist.

 We say $(K(t),0\le t<T)$ is a Loewner chain in $\HH$
w.r.t.\ $\infty$, if each $K(t)$ is a hull in $\HH$ w.r.t.\
$\infty$; $K(0)=\emptyset$; $K(s)\subsetneqq K(t)$ if $s<t$; and for
each fixed $a\in(0,T)$, the extremal length (\cite{Ahl}) of the
curve in $\HH\sem K(t+\eps)$ that disconnect $K(t+\eps)\sem K(t)$
from $\infty$ tends to $0$ as $\eps\to 0^+$, uniformly in
$t\in[0,a]$. If $u(t)$, $0\le t<T$, is a continuous (strictly)
increasing function, and satisfies $u(0)=0$, then
$(K(u^{-1}(t)),0\le t<u(T))$ is also a Loewner chain in $\HH$
w.r.t.\ $\infty$, where $u(T):=\sup u([0,T))$. It is called the
time-change of $(K(t))$ through $u$. Here is a simple example of
Loewner chain. Suppose $\beta(t)$, $0\le t<T$, is a simple curve
with $\beta(0)\in\R$ and $\beta(t)\in\HH$ for $t\in(0,T)$. Let
$K(t)=\beta((0,t])$ for $0\le t<T$. Then $(K(t),0\le t<T)$ is a
Loewner chain in $\HH$ w.r.t.\ $\infty$. It is called the Loewner
chain generated by $\beta$.

If $H_1\subset H_2$ are two hulls in $\HH$ w.r.t.\ $\infty$, let
$H_2/H_1:=\vphi_{H_1}(H_2\sem H_1)$. Then $H_2/H_1$ is also a hull
in $\HH$ w.r.t.\ $\infty$,
$\vphi_{H_2/H_1}=\vphi_{H_2}\circ\vphi_{H_1}^{-1}$, and
$\hcap(H_2/H_1)=\hcap(H_2)-\hcap(H_1)$. If $H_1\subset H_2\subset
H_3$ are three hulls in $\HH$ w.r.t.\ $\infty$, then $H_2/H_1\subset
H_3/H_1$ and $(H_3/H_1)/(H_2/H_1)=H_3/H_2$.

\begin{Proposition} (a) Suppose $K(t)$ and $\vphi(t,\cdot)$, $0\le t<T$,
are chordal Loewner hulls and maps, respectively, driven by $\xi\in
C([0,T))$. Then $(K(t),0\le t<T)$ is a Loewner chain in $\HH$
w.r.t.\ $\infty$, $\vphi_{K(t)}=\vphi(t,\cdot)$, and
$\hcap(K(t))=2t$ for any $0\le t<T$. Moreover, for every
$t\in[0,T)$, \BGE \{\xi(t)\}=\bigcap_{\eps\in(0,T-t)}\lin{K(t+\eps)/
K(t)}.\label{chordal-driving}\EDE (b) Let $(L(s),0\le s<S)$ be a
Loewner chain in $\HH$ w.r.t.\ $\infty$. Let $v(s)=\hcap(L(s))/2$,
$0\le s<S$. Then $v$ is a continuous increasing function with
$u(0)=0$. Let $T=v(S)$ and $K(t)=L(v^{-1}(t))$, $0\le t<T$. Then
$K(t)$, $0\le t<T$, are chordal Loewner hulls driven by some $\xi\in
C([0,T))$. \label{chordal-Loewner-chain}
\end{Proposition}
{\bf Proof.} This is almost the same as Theorem 2.6 in \cite{LSW1}.
$\Box$

\vskip 4mm

Let $B(t)$ be a (standard linear) Brownian motion,
$\kappa\in(0,\infty)$, and $\xi(t)=\sqrt\kappa B(t)$, $0\le
t<\infty$. Let $K(t)$ and $\vphi(t,\cdot)$, $0\le t<\infty$, be the
chordal Loewner hulls and maps, respectively, driven by $\xi$. Then
we call $(K(t))$ the standard chordal $\SLE_\kappa$. From
\cite{RS-basic}, the chordal Loewner trace $\beta(t)$, $0\le
t<\infty$, driven by $\xi$ exists almost surely. Such $\beta$ is
called the standard chordal $\SLE_\kappa$ trace. We have
$\beta(0)=0$ and $\lim_{t\to\infty}\beta(t)=\infty$. If
$\kappa\in(0,4]$, then $\beta$ is a simple curve, $\beta(t)\in\HH$
for $t>0$, and $K(t)=\beta((0,t])$ for $t\ge 0$; if
$\kappa\in(4,\infty)$, then $\beta$ is not a simple curve. If
$\kappa\in[8,\infty)$, then $\beta$ visits every $z\in\lin{\HH}$; if
$\kappa\in (0,8)$, then the Lebesgue measure of the image of $\beta$
in $\C$ is $0$.

Suppose $D\subsetneqq \C$ is a simply connected domain, and $a\ne b$
are two prime ends (\cite{Ahl}) of $D$. Then there is $W$ that maps
$(\HH;0,\infty)$ conformally onto $(D;a,b)$. We call the image of
the standard chordal $\SLE_\kappa$ under $W$ the chordal
$\SLE_\kappa$ in $D$ from $a$ to $b$, which is denoted by
$\SLE_\kappa(D;a\to b)$. Such $W$ is not unique, but the
$\SLE_\kappa(D;a\to b)$ defined through different $W$ have the same
distribution up to a linear time-change because the standard chordal
$\SLE_\kappa$ satisfies the scaling property. The main theorem in
this paper is as follows.

\begin{Theorem} Suppose $\kappa\in(0,4]$, $\beta_1(t)$, $0\le t<\infty$,
 is an $\SLE_\kappa(D;a\to b)$ trace, and $\beta_2(t)$, $0\le t<\infty$,
  is an $\SLE_\kappa(D;b\to a)$ trace. Then the set $\{\beta_1(t):0<
  t<\infty\}$ has the same distribution as $\{\beta_2(t):0<
  t<\infty\}$. \label{Theorem}
\end{Theorem}

\section{Ensemble of Two Chordal Loewner Chains} \label{ensem}
In this section we study the relations of two chordal Loewner chains
that grow together. Some computations were done in \cite{LSW1},
\cite{LSW-8/3}, \cite{Julien}, and other papers. We will give
self-contained arguments for all results in this section. Suppose
$K_j(t)$ and $\vphi_j(t,\cdot)$, $0\le t<S_j$, are chordal Loewner
hulls and maps, respectively driven by $\xi_j\in C([0,S_j))$,
$j=1,2$. Assume that for any $t_1\in[0,S_1)$ and $t_2\in[0,S_2)$,
$\lin{K_1(t_1)}\cap\lin{K_2(t_2)}=\emptyset$, then $K_1(t_1)\cup
K_2(t_2)$ is a hull in $\HH$ w.r.t.\ $\infty$. Fix $j\ne
k\in\{1,2\}$ and $t_0\in[0,S_k)$. For $0\le t<S_j$, let \BGE
K_{j,t_0}(t)=(K_j(t)\cup
K_{k}(t_0))/K_k(t_0)=\vphi_{k}(t_0,K_j(t)).\label{K-j-t}\EDE Since
$\vphi_{k}(t_0,\cdot)$ maps $\HH\sem K_{k}(t_0)$ conformally onto
$\HH$, so from conformal invariance of extremal length,
$(K_{j,t_0}(t),0\le t<S_j)$ is also a Loewner chain in $\HH$ w.r.t.\
$\infty$. Let $v_{j,t_0}(t)=\hcap(K_{j,t_0}(t))/2$ for $0\le t<S_j$,
and $L_{j,t_0}(t)=K_{j,t_0}(v_{j,t_0}^{-1}(t))$ for $0\le
t<S_{j,t_0}:=v_{j,t_0}(S_j)$. From Proposition
\ref{chordal-Loewner-chain}, $L_{j,t_0}(t)$, $0\le t<S_{j,t_0}$, are
chordal Loewner hulls driven by some $\eta_{j,t_0}\in
C([0,S_{j,t_0}))$. Let $\psi_{j,t_0}(t,\cdot)$, $0\le t<S_{j,t_0}$,
denote the corresponding chordal Loewner maps. Let
$\xi_{j,t_0}(t)=\eta_{j,t_0}(v_{j,t_0}(t))$ and
$\vphi_{j,t_0}(t,\cdot)=\psi_{j,t_0}(v_{j,t_0}(t),\cdot)$ for $0\le
t<S_j$. Since $\psi_{j,t_0}(t,\cdot)=\vphi_{L_{j,t_0}(t)}$ for $0\le
t<S_{j,t_0}$, so $\vphi_{j,t_0}(t,\cdot)=\vphi_{K_{j,t_0}(t)}$ for
$0\le t<S_j$. We use $\pa_1$ and $\pa_z$ to denote the partial
derivatives of $\vphi_j(\cdot,\cdot)$ and
$\vphi_{j,t_0}(\cdot,\cdot)$ w.r.t.\ the first (real) and second
(complex) variables, respectively, inside the bracket; and use
$\pa_0$ to denote the partial derivative of
$\vphi_{j,t_0}(\cdot,\cdot)$ w.r.t. the subscript $t_0$.

Fix $j\ne k\in\{1,2\}$, $t\in[0,S_j)$ and $s\in[0,S_{k})$. Since
$\vphi_k(s,\cdot)=\vphi_{K_k(s)}$,
$\vphi_j(t,\cdot)=\vphi_{K_j(t)}$,
$\vphi_{j,s}(t,\cdot)=\vphi_{K_{j,s}(t)}$, and
$\vphi_{k,t}(s,\cdot)=\vphi_{K_{k,t}(s)}$, so from (\ref{K-j-t}),
for any $z\in \HH\sem(K_j(t)\cup K_{k}(s))$, \BGE \vphi_{K_j(t)\cup
K_{k}(s)}(z)= \vphi_{k,t}(s,\vphi_j(t,z))=
\vphi_{j,s}(t,\vphi_{k}(s,z)).\label{circ=1}\EDE Fix
$\eps\in(0,S_j-t)$. Since $K_{j,s}(r)=(K_j(r)\cup
K_{k}(s))/K_{k}(s)$ for $r\in[0,S_j)$, so \BGEN
\frac{L_{j,s}(v_{j,s}(t+\eps))}{L_{j,s}(v_{j,s}(t))}
=\frac{K_{j,s}(t+\eps)}{K_{j,s}(t)}=\frac{K_j(t+\eps)\cup
K_{k}(s)}{K_j(t)\cup K_{k}(s)}\EDEN \BGE =\vphi_{K_j(t)\cup
K_{k}(s)}(K_j(t+\eps)\sem K_j(t))=
\vphi_{k,t}(s,K_j(t+\eps)/K_j(t)).\label{K/K}\EDE From Proposition
\ref{chordal-Loewner-chain} and (\ref{K/K}), we have \BGE
\{\xi_j(t)\}=\cap_{\eps>0}\lin{K_j(t+\eps)/K_j(t)};\quad\mbox{and}\label{xi-j}\EDE
\BGE
\{\xi_{j,s}(t)\}=\{\eta_{j,s}(v_{j,s}(t))\}=\cap_{\eps>0}\lin{L_{j,s}
(v_{j,s}(t+\eps))/L_{j,s}(v_{j,s}(t))}\label{xi-j-t}\EDE \BGE=
\cap_{\eps>0}\lin{(K_j(t+\eps)\cup K_k(s))/(K_j(t)\cup
K_k(s))}.\label{xi-j-t*}\EDE From (\ref{K/K}$\sim$\ref{xi-j-t}), we
have \BGE \xi_{j,s}(t)=\vphi_{k,t}(s,\xi_j(t)).\label{xi=}\EDE From
Proposition \ref{chordal-Loewner-chain} again, we have
$\hcap(K_j(t+\eps)/K_j(t))=2\eps$ and
$$\hcap(L_{j,s}(v_{j,s}(t+\eps))/L_{j,s}(v_{j,s}(t)))=2(v_{j,s}(t+\eps)-v_{j,s}(t)).$$
So from Proposition \ref{hcap} and (\ref{K/K}), we have \BGE
v_{j,s}'(t)=\pa_z\vphi_{k,t}(s,\xi_j(t))^2.\label{v'}\EDE Since
$\vphi_{j,s}(t,z)=\psi_{j,s}(v_{j,s}(t),z)$, so for fixed
$s\in[0,S_k)$, $(t,z)\mapsto \vphi_{j,s}(t,z)$ is $C^{1,a}$
differentiable, where the superscript ``$a$'' means analytic, and
\BGE\pa_1\vphi_{j,s}(t,z)=\frac{2v_{j,s}'(t)}{\psi_{j,s}(v_{j,s}(t),z)
-\eta_{j,s}(v_{j,s}(t))}=\frac{2\pa_z\vphi_{k,t}(s,\xi_j(t))^2}
{\vphi_{j,s}(t,z)-\vphi_{k,t}(s,\xi_j(t))}.\label{chordal*}\EDE From
(\ref{circ=1}), we see that $(s,t,z)\mapsto \vphi_{j,s}(t,z)$ is
$C^{1,1,a}$ differentiable. Differentiate (\ref{chordal*}) using
$\pa_z$, and then divide both sides by $\pa_z\vphi_{j,s}(t,z)$. We
get \BGE
\pa_1\ln(\pa_z\vphi_{j,s}(t,z))=\frac{-2\pa_z\vphi_{k,t}(s,\xi_j(t))^2}
{(\vphi_{j,s}(t,z)-\vphi_{k,t}(s,\xi_j(t)))^2}.\label{chordal**}\EDE
Differentiate (\ref{chordal**}) using $\pa_z$. We get \BGE
\pa_1\Big(\frac{\pa_z^2\vphi_{j,s}(t,z)}{\pa_z\vphi_{j,s}(t,z)}\Big)
=\frac{4\pa_z\vphi_{k,t}(s,\xi_j(t))^2\pa_z\vphi_{j,s}(t,z)}
{(\vphi_{j,s}(t,z)-\vphi_{k,t}(s,\xi_j(t)))^3}.\label{chordal*2}\EDE
Differentiate (\ref{chordal*2}) using $\pa_z$. We get \BGE
\pa_1\pa_z\Big(\frac{\pa_z^2\vphi_{j,s}(t,z)}{\pa_z\vphi_{j,s}(t,z)}\Big)
=\frac{4\pa_z\vphi_{k,t}(s,\xi_j(t))^2\pa_z^2\vphi_{j,s}(t,z)}
{(\vphi_{j,s}(t,z)-\vphi_{k,t}(s,\xi_j(t)))^3}-\frac{12\pa_z\vphi_{k,t}(s,\xi_j(t))^2\pa_z\vphi_{j,s}(t,z)^2}
{(\vphi_{j,s}(t,z)-\vphi_{k,t}(s,\xi_j(t)))^4}.\label{chordal*3}\EDE

\begin{Lemma} For any $j\ne k\in\{0,1\}$, $t\in[0,S_j)$ and $s\in
[0,S_k)$, we have \BGE
\pa_0\vphi_{k,t}(s,\xi_j(t))=-3\pa_z^2\vphi_{k,t}(s,\xi_j(t));\label{-3}\EDE
\BGE
\frac{\pa_0\pa_z\vphi_{k,t}(s,\xi_j(t))}{\pa_z\vphi_{k,t}(s,\xi_j(t))}
=\frac
12\cdot\left(\frac{\pa_z^2\vphi_{k,t}(s,\xi_j(t))}{\pa_z\vphi_{k,t}(s,\xi_j(t))}
\right)^2-\frac
43\cdot\frac{\pa_z^3\vphi_{k,t}(s,\xi_j(t))}{\pa_z\vphi_{k,t}(s,\xi_j(t))}.
\label{1/2-4/3}\EDE \label{-3-1/2-4/3}
\end{Lemma}
{\bf Proof.} Differentiate both sides of the second ``$=$'' in
(\ref{circ=1}) w.r.t.\ $t$, we get
$$\pa_0\vphi_{k,t}(s,\vphi_j(t,z))+\pa_z\vphi_{k,t}(s,\vphi_j(t,z))\pa_1\vphi_j(t,z)
=\pa_1\vphi_{j,s}(t,\vphi_{k}(s,z))$$ for any
$z\in\HH\sem(K_j(t)\cup K_{k}(s))$. So from
(\ref{chordal-equation}), (\ref{circ=1}), and (\ref{chordal*}),
$$\pa_0\vphi_{k,t}(s,\vphi_j(t,z))=\frac{2\pa_z\vphi_{k,t}(s,\xi_j(t))^2}
{\vphi_{k,t}(s,\vphi_{j}(t,z))-\vphi_{k,t}(s,\xi_j(t))}-\frac{2\pa_z\vphi_{k,t}(s,\vphi_j(t,z))}
{\vphi_j(t,z)-\xi_j(t)}$$ for any $z\in\HH\sem(K_j(t)\cup
K_{k}(s))$. Since $\vphi_j(t,\cdot)$ maps $\HH\sem(K_j(t)\cup
K_{k}(s))$ conformally onto $\HH\sem K_{k,t}(s)$, so for any
$w\in\HH\sem K_{k,t}(s)$, \BGE
\pa_0\vphi_{k,t}(s,w)=\frac{2\pa_z\vphi_{k,t}(s,\xi_j(t))^2}
{\vphi_{k,t}(s,w)-\vphi_{k,t}(s,\xi_j(t))}-\frac{2\pa_z\vphi_{k,t}(s,w)}
{w-\xi_j(t)}.\label{pa0-1}\EDE In the above equation, let
$w\to\xi_j(t)$ in $\HH\sem K_{k,t}(s)$. From the Taylor expansion of
$\vphi_{k,t}(s,\cdot)$ at $\xi_j(t)$, we get (\ref{-3}).
Differentiate (\ref{pa0-1}) using $\pa_z$, we get
$$\pa_0\pa_z\vphi_{k,t}(s,w)=-\frac{2\pa_z\vphi_{k,t}(s,\xi_j(t))^2\pa_z\vphi_{k,t}(s,w)}
{(\vphi_{k,t}(s,w)-\vphi_{k,t}(s,\xi_j(t)))^2}-\frac{2\pa_z^2\vphi_{k,t}(s,w)}{w-\xi_j(t)}
+\frac{2\pa_z\vphi_{k,t}(s,w)} {(w-\xi_j(t))^2}.$$ Let
$w\to\xi_j(t)$ in $\HH\sem K_{k,t}(s)$, then we get (\ref{1/2-4/3})
from Taylor expansion. $\Box$

\section{Two-Dimensional Continuous Local Martingale}
\label{Two-Dim}
Let $\kappa\in(0,4]$ and $x_1< x_2\in\R$. Let $X_1(t)$ and $X_2(t)$
be two independent Bessel process of dimension $3-8/\kappa$ started
from $(x_2-x_1)/\sqrt\kappa$. Let $T_j$ denote the first time that
$X_j(t)$ visits $0$, which exists and is finite because
$3-8/\kappa\le 1$. For $j=1,2$, let $Y_j(t)=\sqrt\kappa X_j(t)$,
$0\le t\le T_j$. Then there are two independent Brownian motions
$B_1(t)$ and $B_2(t)$ such that for $j=1,2$ and $0\le t\le T_j$,
$$Y_j(t)=(x_2-x_1)+(-1)^j\sqrt\kappa B_j(t)+\int_0^t\frac{\kappa-4}{Y_j(s)}\,ds.$$
 Fix $j\ne k\in\{1,2\}$. For $0\le t\le
T_j$, let
$$\xi_j(t)=x_j+\sqrt\kappa
B_j(t)+(-1)^j\int_0^t\frac{\kappa-6}{Y_j(s)}\,ds,\quad
p_j(t)=x_k-(-1)^j\int_0^t\frac{2}{Y_j(s)}\,ds.$$ Then
$\xi_j(0)=x_j$, $p_j(0)=x_k$, and $\xi_j(t)-p_j(t)=(-1)^jY_j(t)$,
$0\le t\le T_j$. Thus \BGE d\xi_j(t)=\sqrt\kappa d
B_j(t)+\frac{\kappa-6}{\xi_j(t)-p_j(t)}\,dt,\,\mbox{ and }\, d
p_j(t)=\frac{2dt}{p_j(t)-\xi_j(t)}\label{SDE-xi}\EDE for $0\le t<T$.
Let $K_j(t)$ and $\vphi_j(t,\cdot)$, $0\le t\le T_j$, denote the
chordal Loewner hulls and maps driven by $\xi_j(t)$, $0\le t\le
T_j$. Then $(K_j(t),0\le t<T_j)$ are SLE$(\kappa,\kappa-6)$ process
(\cite{LSW-8/3}) started from $x_j$ with force point at $x_k$; $T_j$
is the first time that $x_k$ is swallowed by $K_j(t)$; and
$\vphi_j(t,x_k)=p_j(t)$, $0\le t< T_j$. It is well known (e.g.\
\cite{Julien}) that after a time-change, $(K_j(t),0\le t<T_j)$ has
the same distribution as a chordal $\SLE_\kappa(\HH;x_j\to x_k)$.
Since $\kappa\le 4$, so there is a crosscut $\beta_j(t)$, $0\le t\le
T_j$, in $\HH$ from $x_j$ to $x_k$, such that
$K_j(t)=\beta_j((0,t])$ for $0\le t< T_j$ (\cite{RS-basic}). Here a
crosscut in $\HH$ from $a\in\R$ to $b\in\R$ is simple curve
$\beta(t)$, $0\le t\le T$, that satisfies $\beta(0)=a$,
$\beta(T)=b$, and $\beta(t)\in\HH$ for $0<t<T$.

For $j=1,2$, let $(\F^j_t)$ denote the filtration generated by
$(B_j(t))$. Then $(\xi_j)$ is $(\F^j_t)$-adapted, and $T_j$ is an
$(\F^j_t)$-stopping time. Let
$${\cal D}=\{(t_1,t_2)\in[0,T_1)\times[0,T_2):\lin{K_1(t_1)}\cap\lin{K_2(t_2)}=\emptyset\}.$$
For $0\le t_k<T_k$, let $T_j(t_k)\in(0,T_j]$ be the maximal such
that $\lin{K_j(t)}\cap\lin{K_k(t_k)}\ne\emptyset$ for $0\le
t<T_j(t_k)$. Now we use the notations in the last section. Let
$(t_1,t_2)\in {\cal D}$. Since
$\vphi_{k,t_j}(t_k,\cdot)=\vphi_{K_{k,t_j}(t_k)}$, so
$\vphi_{k,t_j}(t_k,\cdot)$ maps $\HH\sem K_{k,t_j}(t_k)$ conformally
onto $\HH$. By Schwarz reflection principle,
$\vphi_{k,t_j}(t_k,\cdot)$ extends conformally to
$\Sigma_{K_{k,t_j}(t_k)}$, where for a hull $H$ in $\HH$ w.r.t.\
$\infty$, $\Sigma_H=\C\sem(H\cup \{\lin{z}:z\in H\}\cup
[\inf(\lin{H}\cap\R),\sup(\lin{H}\cap\R)])$ (c.f.\ \cite{LERW}). For
$j\ne k\in\{0,1\}$, and $h\in\Z_{\ge 0}$, let
$A_{j,h}(t_1,t_2)=\pa_z^h\vphi_{k,t_j}(t_k,\xi_j(t_j))$. The
definition makes sense since $\xi_j(t_j)\in
\Sigma_{K_{k,t_j}(t_k)}$. Moreover, we have $A_{j,h}\in\R$ for any
$h\ge 0$ since $\vphi_{k,t_j}(t_k,\cdot)$ is real valued on a real
open interval containing $\xi_j(t_j)$.  From (\ref{circ=1}), we see
that $A_{j,0}(t_1,t_2)=\vphi_{K_1(t_1)\cup K_2(t_2)}(\beta_j(t_j))$,
$j=1,2$. Since $K_1(t_1)$ lies to the left of $K_2(t_2)$, so
$A_{1,0}(t_1,t_2)<A_{2,0}(t_1,t_2)$. Since
$\vphi_{k,t_j}(t_k,\cdot)$ maps a part of the upper half plane to
the upper half plane, so $A_{j,1}(t_1,t_2)>0$, $j=1,2$. For
$(t_1,t_2)\in {\cal D}$, define
$E(t_1,t_2)=A_{2,0}(t_1,t_2)-A_{1,0}(t_1,t_2)>0$, \BGE
N(t_1,t_2)=\frac{A_{1,1}(t_1,t_2)A_{2,1}(t_1,t_2)}{E(t_1,t_2)^2}
=\frac{A_{1,1}(t_1,t_2)A_{2,1}(t_1,t_2)}
{(A_{2,0}(t_1,t_2)-A_{1,0}(t_1,t_2))^2}>0,\label{N}\EDE and \BGE
M(t_1,t_2)=\Big(\frac{N(t_1,t_2)N(0,0)}{N(t_1,0)N(0,t_2)}\Big)^\alpha\exp\Big(-\lambda\int_0^{t_1}\!\!\int_0^{t_2}
2N(s_1,s_2)^2\,ds_2\,ds_1\Big)>0,\label{M-def}\EDE where \BGE
\alpha=\alpha(\kappa)=\frac{6-\kappa}{2\kappa},\quad
\lambda=\lambda(\kappa)=\frac{(8-3\kappa)(6-\kappa)}{2\kappa}.\label{alpha-lambda}\EDE
Note that $M(t_1,0)=M(0,t_2)=1$ for any $0\le t_1<T_1$ and $0\le
t_2<T_2$.

\vskip 3mm

\no {\bf Remark.} If $\kappa<8/3$, i.e., $\lambda>0$, then
$$\exp\Big(-\lambda\int_0^{t_1}\!\!\int_0^{t_2}
2N(s_1,s_2)^2\,ds_2ds_1\Big)$$ is the probability that in a loop
soup (\cite{Loop-soup}) in $\HH$ with intensity $\lambda$, there is
no loop that intersects both $K_1(t_1)$ and $K_2(t_2)$.

\begin{Theorem} (i) For any fixed $(\F^2_t)$-stopping time $\bar t_2$ with $\bar t_2<T_2$,
$(M(t_1,\bar t_2),0\le t_1<T_1(\bar t_2))$ is a continuous
$(\F^1_{t_1}\times\F^2_{\bar t_2})_{t_1\ge 0}$-local martingale, and
\BGE \frac{\pa_1 M}{M}\Big|_{(t_1,\bar t_2)}=\Big(3-\frac \kappa
2\Big)\Big(\Big(\frac{A_{1,2}}
{A_{1,1}}+\frac{2A_{1,1}}{A_{2,0}-A_{1,0}}\Big)\Big|_{(t_1,\bar
t_2)}-\frac{2}{p_1(t_1)-\xi_1(t_1)}\Big)\frac{\pa
B_1(t_1)}{\sqrt\kappa}.\label{t1}\EDE (ii) For any fixed
$(\F^1_t)$-stopping time $\bar t_1$ with $\bar t_1<T_1$, $(M(\bar
t_1,t_2),0\le t_2<T_2(\bar t_1))$ is a continuous $(\F^1_{\bar
t_1}\times\F^2_{t_2})_{t_2\ge 0}$-local martingale, and \BGE
\frac{\pa_2 M}{M}\Big|_{(\bar t_1,t_2)}=\Big(3-\frac \kappa
2\Big)\Big(\Big(\frac{A_{2,2}}
{A_{2,1}}+\frac{2A_{2,1}}{A_{1,0}-A_{2,0}}\Big)\Big|_{(\bar t_1,
t_2)}-\frac{2}{p_2(t_2)-\xi_2(t_2)}\Big)\frac{\pa
B_2(t_2)}{\sqrt\kappa}.\label{t2}\EDE \label{M}
\end{Theorem} \vskip -4mm
{\bf Proof.} Since $\vphi_{1,t_1}(0,\cdot)=\id_\HH$,
$\vphi_{1,0}(t_1,\cdot)=\vphi_1(t_1,\cdot)$, and $\xi_2(0)=x_2$, so
$$A_{1,0}(t_1,0)=\vphi_{2,t_1}(0,\xi_1(t_1))=\xi_1(t_1),\quad
A_{1,1}(t_1,0)=1;$$
$$A_{2,0}(t_1,0)=\vphi_{1,0}(t_1,\xi_2(0))=\vphi_1(t_1,x_2)=p_1(t_1),\quad
A_{2,1}(t_1,0)=\pa_z\vphi_1(t_1,x_2).$$ Thus
$N(t_1,0)=\pa_z\vphi_1(t_1,x_2)/(p_1(t_1)-\xi_1(t_1))^2$. From
chordal Loewner equation, we get
$$\pa_{t_1}\pa_z\vphi_1(t_1,x_2)=\frac{-2\pa_z\vphi_1(t_1,x_2)}{(\vphi_1(t_1,x_2)-\xi_1(t_1))^2}
=\frac{-2\pa_z\vphi_1(t_1,x_2)}{(p_1(t_1)-\xi_1(t_1))^2}.$$ From
(\ref{SDE-xi}), we get
$$\pa_{t_1}(p_1(t_1)-\xi_1(t_1))=-\pa \xi_1(t_1)+\frac{2\pa
t_1}{p_1(t_1)-\xi_1(t_1)}.$$ From the above two formulas and Ito's
formula, we get \BGE \pa_1 N(t_1,0)^\alpha/(\alpha
N(t_1,0)^\alpha)=2\pa
\xi_1(t_1)/(p_1(t_1)-\xi_1(t_1)).\label{N-0}\EDE

Now fix an $(\F^2_t)$-stopping time $\bar t_2$ with $\bar t_2<T_2$.
Then we get a filtration $(\F^1_{t}\times\F^2_{\bar t_2})_{t\ge 0}$.
Since $B_1(t)$ and $B_2(t)$ are independent, so $B_1(t)$ is an
$(\F^1_{t}\times\F^2_{\bar t_2})_{t\ge 0}$-Browinian motion. Then
$T_1(\bar t_2)$ is an $(\F^1_{t}\times\F^2_{\bar t_2})_{t\ge
0}$-stopping time, $A_{j,h}(t,\bar t_2)$, $j=1,2$, $E(t,\bar t_2)$,
$N(t,\bar t_2)$, and $M(t,\bar t_2)$ are defined for
$t\in[0,T_1(\bar t_2))$. From chordal Loewner equation and
(\ref{circ=1}), $\vphi_{1,\bar t_2}(t,\cdot))$ and $\vphi_{2,t}(\bar
t_2,\cdot)$, $0\le t<T_1(\bar t_2)$, are $(\F^1_{t}\times\F^2_{\bar
t_2})_{t\ge 0}$-adapted. Since $A_{1,h}(t,\bar
t_2)=\pa_z^h\vphi_{2,t}(\bar t_2,\xi_1(t))$,  so from Ito's formula,
$(A_{1,h}(t_1,\bar t_2),0\le t_1<T_1(\bar t_2))$ satisfies the
$(\F^1_{t}\times\F^2_{\bar t_2})_{t\ge 0}$-adapted SDE: \BGE \pa_1
A_{1,h}(t_1,\bar t_2)= A_{1,h+1}(t_1,\bar t_2)\pa
\xi_1(t_1)+\Big(\pa_0\pa_z^h \vphi_{2,t_1}(\bar t_2,\xi_1(t_1))
+\frac\kappa 2 A_{1,h+2}(t,\bar t_2)\Big)\,\pa t_1.\label{SDE-A}\EDE
From (\ref{chordal*}) and (\ref{chordal**}), we have \BGE \pa_1
A_{2,0}(t_1,t_2)=\frac{2 A_{1,1}(t_1,t_2)^2}{E(t_1,t_2)}\,\pa
t_1,\quad \frac{\pa_1 A_{2,1}(t_1,t_2)}{A_{2,1}(t_1,t_2)}=-\frac{2
A_{1,1}(t_1,t_2)^2}{E(t_1,t_2)^2}\,\pa t_1.\label{pa-A_2}\EDE
 From (\ref{SDE-A}), (\ref{pa-A_2}), and  Lemma
\ref{-3-1/2-4/3}, we have \BGE \pa_1 A_{1,0}= A_{1,1} \pa
\xi_1(t_1)+(\frac\kappa 2-3)A_{1,2}\pa
t_1,\quad\mbox{and}\label{pa-A_1-0}\EDE \BGE \frac{\pa_1
A_{1,1}}{A_{1,1}}= \frac{A_{1,2}} {A_{1,1}}\,\pa
\xi_1(t_1)+\Big(\frac 12\cdot\Big(\frac{A_{1,2}}
{A_{1,1}}\Big)^2+\Big(\frac\kappa 2-\frac 43\Big)\cdot\frac{A_{1,3}}
{A_{1,1}}\Big)\,\pa t_1,\label{pa-A_1-1}\EDE where ``$(t_1,\bar
t_2)$'' are omitted. Since $E=A_{2,0}-A_{1,0}$, from (\ref{pa-A_2})
and (\ref{pa-A_1-0}), we have \BGE \pa_1 E=-A_{1,1}\pa
\xi_1(t_1)+(\frac{2A_{1,1}^2}{E}+(3-\frac\kappa 2)A_{1,2})\,\pa
t_1.\label{pa-D}\EDE  Let $C_h=A_{1,h}$ for $h\in\Z_{\ge 0}$. From
(\ref{pa-A_2}$\sim$\ref{pa-D}) and Ito's formula, we have \BGE
\frac{\pa_1 N^\alpha}{\alpha
N^\alpha}=\Big(\frac{C_2}{C_1}+\frac{2C_1}{E}\Big)\,\pa \xi_1(t_1)
+(8-3\kappa)\Big(\frac 14\cdot\frac{C_2^2}{C_1^2}-\frac
16\cdot\frac{C_3}{C_1}\Big)\,\pa t_1.\label{pa-N-lam}\EDE The above
SDE is $(\F^1_{t}\times\F^2_{\bar t_2})_{t\ge 0}$-adapted. Now
(\ref{N-0}) is also an $(\F^1_{t}\times\F^2_{\bar t_2})_{t\ge
0}$-adapted SDE since $B_1(t)$ is an $(\F^1_{t}\times\F^2_{\bar
t_2})_{t\ge 0}$-Brownian motion. Thus from (\ref{SDE-xi}),
(\ref{N-0}), (\ref{pa-N-lam}), and Ito's formula, we have $$
\frac{\pa_1(N(t_1,\bar t_2)/N(t_1,0))^\alpha} {\alpha(N(t_1,\bar
t_2)/N(t_1,0))^\alpha}=\Big(\frac{C_2(t_1,\bar t_2)}{C_1(t_1,\bar
t_2)}+\frac{2C_1(t_1,\bar t_2)}{E(t_1,\bar t_2)}-\frac
2{p_1(t_1)-\xi_1(t_1)}\Big)\,\sqrt\kappa\pa B_1(t_1)$$ \BGE
+(8-3\kappa)\Big(\frac 14\cdot\frac{C_2(t_1,\bar
t_2)^2}{C_1(t_1,\bar t_2)^2}-\frac 16\cdot\frac{C_3(t_1,\bar
t_2)}{C_1(t_1,\bar t_2)}\Big)\,\pa t_1.\label{pa-N-0}\EDE

Since $C_j(t_1,t_2)=\pa_z^j\vphi_{2,t_1}(t_2,\xi_1(t_1))$, so $\pa_2
C_j(t_1,t_2)=\pa_1\pa_z^j\vphi_{2,t_1}(t_2,\xi_1(t_1))$, and
$$\Big(\frac 14\cdot\frac{C_2^2}{C_1^2}-\frac
1 6\cdot\frac{C_3}{C_1}\Big)\Big|_{(t_1,t_2)}
=\frac{1}{12}({\pa_z^2}/
{\pa_z})\vphi_{2,t_1}(t_2,\xi_1(t_1))^2-\frac 1 6\pa_z({\pa_z^2}/
{\pa_z})\vphi_{2,t_1}(t_2,\xi_1(t_1)).$$ From (\ref{chordal*2}) and
(\ref{chordal*3}), we have
$$\frac{\pa}{\pa t_2}\Big[({\pa_z^2}/
{\pa_z})\vphi_{2,t_1}(t_2,\xi_1(t_1))^2\Big]=\frac{8A_{2,1}^2C_2}{E^3}\Big|_{(t_1,t_2)},$$
$$\frac{\pa}{\pa t_2}\Big[\pa_z({\pa_z^2}/
{\pa_z})\vphi_{2,t_1}(t_2,\xi_1(t_1))\Big]=\Big(\frac{4A_{2,1}^2C_2}{E^3}-
\frac{12A_{2,1}^2C_1^2}{E^4}\Big)\Big|_{(t_1,t_2)}.$$ From the above
three formulas, we get
$$\pa_2\Big(\frac 14\cdot\frac{C_2^2}{C_1^2}-\frac
16\cdot\frac{C_3}{C_1}\Big)\Big|_{(t_1,t_2)}=\frac{2A_{2,1}^2C_1^2}
{E^4}\Big|_{(t_1,t_2)}=2 N(t_1,t_2)^2.$$ Since
$\vphi_{2,t_1}(0,\cdot)=\id_\HH$, so
$\pa_z^j\vphi_{2,t_1}(0,\cdot)=0$ for $j\ge 2$. Thus
$C_2(t_1,0)=C_3(t_1,0)=0$. So \BGE \frac
14\cdot\frac{C_2(t_1,t_2)^2}{C_1(t_1,t_2)^2}-\frac
16\cdot\frac{C_3(t_1,t_2)}{C_1(t_1,t_2)}=\int_0^{t_2}
2N(t_1,s_2)^2ds_2.\label{int}\EDE Then (\ref{t1}) follows from
(\ref{M-def}$\sim$\ref{alpha-lambda}) and
(\ref{pa-N-0}$\sim$\ref{int}); (\ref{t2}) follows from the symmetry.
$\Box$

\vskip 3mm

Now we make some improvement over the above theorem. Let $\bar t_2$
be an $(\F^2_t)$-stopping time with $\bar t_2<T_2$. Suppose $R$ is
an $(\F^1_{t}\times\F^2_{\bar t_2})_{t\ge 0}$-stopping time with $R<
T_1(\bar t_2)$. Let $\F_{R,\bar t_2}$ denote the $\sigma$-algebra
obtained from the filtration $(\F^1_t\times\F^2_{\bar t_2})_{t\ge
0}$ and its stopping time $R$, i.e., ${\cal E}\in \F_{R,\bar t_2}$
iff for any $t\ge 0$, ${\cal E}\cap\{R\le t\}\in
\F^1_t\times\F^2_{\bar t_2}$. For every $t\ge 0$, $R+t$ is also an
$(\F^1_{t}\times\F^2_{\bar t_2})_{t\ge 0}$-stopping time. So we have
a filtration $(\F_{R+t,\bar t_2})_{t\ge 0}$. Since $(\xi_1(t))$ and
$(p_1(t))$ are $(\F^1_{t}\times\F^2_{\bar t_2})_{t\ge 0}$-adapted,
so $(\xi_1(R+t),t\ge 0)$, $(p_1(R+t),t\ge 0)$,
$(\vphi_1(R+t,\cdot),t\ge 0))$ and $(K_1(R+t),t\ge 0)$ are
$(\F_{R+t,\bar t_2})_{t\ge 0}$-adapted. Suppose $I\in[0,\bar t_2]$
is $\F_{R,\bar t_2}$-measurable. From $I\le \bar t_2$ we have
$T_1(I)\ge T_1(\bar t_2)>R$. Then $\vphi_{1,I}(R+t,\cdot)$ and
$\vphi_{2,R+t}(I,\cdot)$ are defined for $0\le t<T_1(I)-R$.

\begin{Lemma} $T_1(I)-R$ is an $(\F_{R+t,\bar t_2})_{t\ge
0}$-stopping time, $(\vphi_{1,I}(R+t,\cdot),0\le t<T_1(I)-R)$ and
$(\vphi_{2,R+t}(I,\cdot),0\le t<T_1(I)-R)$ are $(\F_{R+t,\bar
t_2})_{t\ge 0}$-adapted. \label{measurable}
\end{Lemma}
{\bf Proof.} Since $T_1(I)-R>t$ iff $K_1(R+t)\cap K_2(I)=\emptyset$,
and that $(\vphi_1(R+t,\cdot))$, and $(K_1(R+t))$ are $\F_{R+t,\bar
t_2}$-adapted, so from (\ref{circ=1}), we suffice to show that
$\vphi_{2}(I,\cdot)$ is $\F_{R,\bar t_2}$-measurable. Fix $n\in\N$.
Let $I_n=\lfloor nI\rfloor/n$. For $m\in\N\cup\{0\}$, let ${\cal
E}_n(m)=\{m/n\le I_n<(m+1)/n\}$. Then ${\cal E}_n(m)$ is $\F_{R,\bar
t_2}$-measurable, and $I_n=m/n$ on ${\cal E}_n(m)$. Since $m/n\le
\bar t_2$ and $I_n=m/n$ on ${\cal E}_n(m)$, so $I_n$ agrees with
$(m/n)\wedge \bar t_2$ on ${\cal E}_n(m)$. Now $(m/n)\wedge \bar
t_2$ is an $(\F^2_t)$-stopping time, and $\F^2_{(m/n)\wedge \bar
t_2}\subset\F^2_{\bar t_2}\subset \F_{R,\bar t_2} $. So
$\vphi_2((m/n)\wedge \bar t_2,\cdot)$ is $\F_{R,\bar
t_2}$-measurable. Since $\vphi_2(I_n,\cdot)=\vphi_2((m/n)\wedge \bar
t_2,\cdot)$ on ${\cal E}_n(m)$, and  ${\cal E}_n(m)$ is $\F_{R,\bar
t_2}$-measurable for each $m\in\N\cup\{0\}$, so $\vphi_2(I_n,\cdot)$
is $\F_{R,\bar t_2}$-measurable. Since
$\vphi_2(I_n,\cdot)\to\vphi_2(I,\cdot)$ as $n\to\infty$, so
$\vphi_2(I,\cdot)$ is also $\F_{R,\bar t_2}$-measurable. Then we are
done. $\Box$

\vskip 3mm

Let $B_1^R(t)=B_1(R+t)-B_1(R)$, $0\le t<\infty$. Since $B_1(t)$ is
an $(\F^1_{t}\times\F^2_{\bar t_2})_{t\ge 0}$-Brownian motion, so
$B_1^R(t)$ is an $(\F_{R+t,\bar t_2})_{t\ge 0}$-Brownian motion.
Then $(\xi_1(R+t))$ satisfies the $(\F_{R+t,\bar t_2})_{t\ge
0}$-adapted SDE:
$$d\xi_1(R+t)= \sqrt\kappa
dB^R_1(t)+\frac{\kappa-6}{\xi_1(R+t)-p_1(R+t)}\,dt.$$ The SDEs in
the proof of Theorem \ref{M} still hold if  $t_1$ is replaced by
$R+t$,  $\bar t_2$ is replaced by $I$, and $B_1(t_1)$ is replaced by
$B_1^R(t_1)$. The difference is that the SDEs now are all
$(\F_{R+t,\bar t_2})_{t\ge 0}$-adapted. So we have the following
theorem.

\begin{Theorem} (i) Suppose $\bar t_2$ is an $(\F^2_t)$-stopping
time with $\bar t_2<T_2$. Suppose $R$ is an $(\F^1_{t}\times
\F^2_{\bar t_2})_{t\ge 0}$-stopping time with $R<T_1(\bar t_2)$. Let
$I\in[0,\bar t_2]$ be $\F_{R,\bar t_2}$-measurable. Then $(M(R+t,I),
0\le t<T_1(I)-R)$ is a continuous $(\F_{R+t,\bar t_2})_{t\ge
0}$-local martingale.\\
(ii) Suppose $\bar t_1$ is an $(\F^1_t)$-stopping time with $\bar
t_1<T_1$. Suppose $I$ is an $(\F^1_{\bar t_1}\times \F^2_{t})_{t\ge
0}$-stopping time with $I<T_2(\bar t_1)$. Let $R\in[0,\bar t_1]$ be
$\F_{\bar t_1,I}$-measurable. Then $(M(R,I+t), 0\le t<T_2(R)-I)$ is
a continuous $(\F_{\bar t_1,I+t})_{t\ge 0}$-local martingale.
 \label{M-I}
\end{Theorem}
{\bf Proof.} (i) follows from the above argument. (ii) follows from
the symmetry. $\Box$

\section{Boundedness} \label{sec-bound}
We now use the notations and results in Section 5.2 of \cite{LERW}.
Let $H$ be a nonempty hull in $\HH$ w.r.t.\ $\infty$. Then
$a_H=\inf\{\lin{H}\cap\R\}$, $b_H=\sup\{\lin{H}\cap\R\}$,
$\Sigma_H=\C\sem(H\cup \{\lin{z}:z\in H\}\cup [a_H,b_H])$, and
${\cal H}(H)$ is the set of hulls in $\HH$ w.r.t.\ $\infty$ that are
contained in $H$. From Lemma 5.4 in \cite{LERW}, any sequence
$(K_n)$ in ${\cal H}(H)$ contains a subsequence $(L_n)$ such that
$\vphi_{L_n}\luto \vphi_K$ (converges locally uniformly) in
$\Sigma_H$ for some $K\in{\cal H}(H)$. We now make some improvement
over this result. Let $Q_H=\lin{H}\cap\R$. Then $Q_H$ is a closed
subset of $[a_H,b_H]$. Let
$$\Sigma_H^*=\Sigma_H\cup([a_H,b_H]\sem Q_H)=\C\sem(H\cup\{\lin{z}:z\in H\}\cup
Q_H),$$ which may strictly contains $\Sigma_H$. For any $K\in{\cal
H}(H)$, $\vphi_K$ extends conformally to $\Sigma_H^*$ by Schwarz
reflection principle, and $\vphi_K'(x)>0$ for any $x\in\R\sem Q_H$
from (5.1) in \cite{LERW}, so $\vphi_K$ preserves the order on
$\R\sem Q_H$.

\begin{Lemma} Suppose $(K_n)$ is a sequence in ${\cal H}(H)$. Then
it contains some subsequence $(L_n)$ such that $\vphi_{L_n}\luto
\vphi_K$ in $\Sigma_H^*$ for some $K\in{\cal H}(H)$.
 \label{comp}
\end{Lemma}
{\bf Proof.} From the argument after Corollary 5.1 in \cite{LERW},
there is $M_H>0$ such that $|\vphi_K(z)-z|\le M_H$ for any
$K\in{\cal H}(H)$ and $z\in \Sigma_H$. After the extension, we have
$|\vphi_K(z)-z|\le M_H$ for any $K\in{\cal H}(H)$ and
$z\in\Sigma_H^*$. So $\{\vphi_{K_n}(z)-z:n\in\N\}$ is a normal
family in $\Sigma_H^*$. Then $(K_n)$ contains a subsequence $(L_n)$
such that $\vphi_{L_n}(z)-z\luto f(z)$ in $\Sigma_H^*$ for some $f$
that is analytic in $\Sigma_H^*$. So $\vphi_{L_n}\luto g$ in
$\Sigma_H^*$, where $g(z):=z+f(z)$ is analytic in $\Sigma_H^*$. From
Lemma 5.4 in \cite{LERW}, we may assume that
$\vphi_{L_n}\luto\vphi_K$ in $\Sigma_H$ for some $K\in{\cal H}(H)$.
Thus $g=\vphi_K$ in $\Sigma_H$. Since they are both analytic in
$\Sigma_H^*$, so $g=\vphi_K$ in $\Sigma_H^*$. Thus $\vphi_{L_n}\luto
\vphi_K$ in $\Sigma_H^*$. $\Box$

\begin{Lemma} If $y_1<y_2<a_H$ or $y_1>y_2>b_H$, then
$\vphi_H'(y_1)> \vphi_H'(y_2)$. \label{order}
\end{Lemma}
{\bf Proof.} This follows from differentiating (5.1) in \cite{LERW}
for $z\in\R\sem[c_H,d_H]$, and the facts that $\vphi_H$ is
increasing on $(-\infty,a_H)$ and $(b_H,\infty)$, and maps them to
$(-\infty,c_H)$ and $(d_H,\infty)$, respectively. $\Box$ \vskip 3mm

Let $\HP$ denote the set of $(H_1,H_2)$ such that $H_j$ is a hull in
$\HH$ w.r.t.\ $\infty$ that contains some neighborhood of $x_j$ in
$\HH$, $j=1,2$, and $\lin{H_1}\cap\lin{H_2}=\emptyset$. Let
$(H_1,H_2)\in\HP$. Then $b_{H_1}<a_{H_2}$, $H_1\cup H_2$ is a
 hull in $\HH$ w.r.t.\ $\infty$, and
$Q_{H_1\cup H_2}=Q_{H_1}\cup
Q_{H_2}\subset[a_{H_1},b_{H_1}]\cup[a_{H_2},b_{H_2}]$. Let
$T_j(H_j)$ be the first time that $\lin{K_j(t)}\cap\lin{\HH\sem
H_j}\ne\emptyset$, $j=1,2$. Then $T_j(H_j)$ is an
$(\F^j_t)$-stopping time, $0<T_j(H_j)<T_j$, and $K_j(t)\subset H_j$
for $0\le t\le T_j({H_j})$. Thus \BGE
T_j(H_j)=\hcap(K_j(T_j(H_j)))/2\le \hcap(H_j)/2.\label{TjH}\EDE

\begin{Theorem} For any $(H_1,H_2)\in\HP$, there are $C_2>C_1>0$
depending only on $H_1$ and $H_2$ such that $M(t_1,t_2)\in[C_1,C_2]$
for any $(t_1,t_2)\in[0,T_1(H_1)]\times [0,T_2(H_2)]$. \label{bound}
\end{Theorem}
{\bf Proof.} Let $(H_1,H_2)\in\HP$ and $H=H_1\cup H_2$. Throughout
this proof, we use $C_n$, $n\in\N$, to denote some positive constant
that depends only on $H_1$ and $H_2$. From (\ref{M-def}) and
(\ref{TjH}), we suffice to show that for some $C_4>C_3>0$,
$N(t_1,t_2)\in [C_3,C_4]$ for $(t_1,t_2)\in[0,T_1(H_1)]\times
[0,T_2(H_2)]$. Fix $(t_1,t_2)\in[0,T_1(H_1)]\times [0,T_2(H_2)]$.
First suppose $t_1,t_2>0$. Fix $j\ne k\in\{1,2\}$. For any
$s_j\in[0,t_j)$, from (\ref{xi-j}) we have
$\xi_j(s_j)\in\lin{K_j(t_j)/K_j(s_j)}$, so
$$\xi_j(s_j)\in[a_{K_j(t_j)/K_j(s_j)},b_{K_j(t_j)/K_j(s_j)}]
\subset[c_{K_j(t_j)/K_j(s_j)},d_{K_j(t_j)/K_j(s_j)}]\subset
[c_{K_j(t_j)},d_{K_j(t_j)}],$$ where the second and third inclusions
follow from Lemma 5.2 and Lemma 5.3 in \cite{LERW}. Let $s_j\to
t_j$. We get $\xi_j(t_j)\in[c_{K_j(t_j)},d_{K_j(t_j)}]$. For
$s_j\in[0,t_j)$, from (\ref{xi-j-t*}) and (\ref{xi=}),
$$A_{j,0}(s_j,t_k)=\vphi_{k,s_j}(t_k,\xi_{j}(s_j))\in\lin{
{(K_j(t_j)\cup K_k(t_k))}/{(K_j(s_j)\cup K_k(t_k))}},$$ which
implies that $A_{j,0}(s_j,t_k)\in[c_{K_j(t_j)\cup K_k(t_k)},
d_{K_j(t_j)\cup K_k(t_k)}]\subset[c_H,d_H]$. Let $s_j\to t_j$. We
get $A_{j,0}(t_j,t_k)\in[c_H,d_H]$. This also holds for
$A_{k,0}(t_j,t_k)$. Thus \BGE
|E(t_j,t_k)|=|A_{j,0}(t_j,t_k)-A_{k,0}(t_j,t_k)|\le
d_H-c_H.\label{Du}\EDE

Fix $q_1,q_2,r_1,r_2\in\R$ with $r_1<a_{H_1}\le
b_{H_1}<q_1<q_2<a_{H_2}\le b_{H_2}<r_2$. From Lemma \ref{comp},
there are $C_6>C_5>0$ such that for $x=q_1,q_2,r_1,r_2$,
$\pa_z\vphi_{K_1(t_1)\cup K_2(t_2)}(x)$, $\pa_z\vphi_1(t_1,x)$, and
$\pa_z\vphi_2(t_2,x)$ all lie in $[C_5,C_6]$. Fix $j\ne
k\in\{1,2\}$. From (\ref{circ=1}) there are $C_8>C_7>0$ such that
for $x=q_j,r_j$,
$\pa_z\vphi_{k,t_j}(t_k,\vphi_j(t_j,x))\in[C_7,C_8]$. Since
$[a_{K_j(t_j)},b_{K_j(t_j)}]\subset[a_{H_j},b_{H_j}]$, so $r_j$ is
disconnected from $q_j$ in $\R$ by $[a_{K_j(t_j)},b_{K_j(t_j)}]$.
Since $\vphi_j(t_j,\cdot)=\vphi_{K_j(t_j)}$, so
 $\vphi_j(t_j,r_j)$ is disconnected from $\vphi_j(t_j,q_j)$ in $\R$ by
$[c_{K_j(t_j)},d_{K_j(t_j)}]$. Since
$\xi_j(t_j)\in[c_{K_j(t_j)},d_{K_j(t_j)}]$, so $\xi_j(t_j)$ lies
between $\vphi_j(t_j,r_j)$ and $\vphi_j(t_j,q_j)$. Since $r_j$ and
$q_j$ lie on the same side of $K_k(t_k)$, so $\vphi_j(t_j,r_j)$,
$\xi_j(t_j)$, and $\vphi_j(t_j,q_j)$ lie on the same side of
$\vphi_j(t_j,K_k(t_k))=K_{k,t_j}(t_k)$. Since
$\vphi_{k,t_j}(t_k,\cdot)=\vphi_{K_{k,t_j}(t_k)}$, so from Lemma
\ref{order}, $\pa_z \vphi_{k,t_j}(t_k,\xi_j(t_j))$ lies between
$\pa_z \vphi_{k,t_j}(t_k,\vphi_j(t_j,r_j))$ and $\pa_z
\vphi_{k,t_j}(t_k,\vphi_j(t_j,q_j))$. Thus \BGE
A_{j,1}(t_j,t_k)=\pa_z
\vphi_{k,t_j}(t_k,\xi_j(t_j))\in[C_7,C_8].\label{A1}\EDE From
(\ref{circ=1}) and the above argument, we see that
$A_{j,0}(t_j,t_k)=\vphi_{k,t_j}(t_k,\xi_j(t_j))$ lies between
$\vphi_{K_j(t_j)\cup K_k(t_k)}(r_j)$ and $\vphi_{K_j(t_j)\cup
K_k(t_k)}(q_j)$ for $j=1,2$. Since $r_1<q_1<q_2<r_2$, so
$$\vphi_{K_1(t_2)\cup
K_2(t_2)}(r_1)<\vphi_{K_1(t_1)\cup
K_2(t_2)}(q_1)<\vphi_{K_1(t_1)\cup
K_2(t_2)}(q_2)<\vphi_{K_1(t_1)\cup K_2(t_2)}(r_2);$$
 From Lemma \ref{comp}, there
is $C_9>0$ such that $\pa_z\vphi_{K_1(t_1)\cup K_2(t_2)}(x)\ge C_9$
for $x\in[q_1,q_2]$. So \BGE |E(t_1,t_2)|\ge \vphi_{K_1(t_1)\cup
K_2(t_2)}(q_2) -\vphi_{K_1(t_1)\cup K_2(t_2)}(q_1)\ge
C_9(q_2-q_1).\label{Dl}\EDE From (\ref{Du}), (\ref{A1}), and
(\ref{Dl}), we have $C_4>C_3>0$ such that $N(t_1,t_2)\in [C_3,C_4]$
for $(t_1,t_2)\in(0,T_1(H_1)]\times(0,T_2(H_2)]$. By letting $t_1$
or $t_2$ tend to $0$, we obtain the above inequality in the case
$t_1$ or $t_2$ equals to $0$. So we are done. $\Box$

\vskip 3mm

Now we explain the meaning of $M(t_1,t_2)$. Fix $(H_1,H_2)\in\HP$.
Let $\mu$ denote the joint distribution of $(\xi_1(t):0\le t\le
T_1)$ and $(\xi_2(t):0\le t\le T_2)$. From Theorem \ref{M} and
Theorem \ref{bound}, we have $\int
M(T_1(H_1),T_2(H_2))d\mu=\EE[M(T_1(H_1),T_2(H_2))]=M(0,0)=1$. Note
that $M(T_1(H_1),T_2(H_2))>0$. Suppose $\nu$ is a measure on
$\F^1_{T_1(H_1)}\times\F^2_{T_2(H_2)}$ such that
$d\nu/d\mu=M(T_1(H_1),T_2(H_2))$. Then $\nu$ is a probability
measure. Now suppose the joint distribution of $(\xi_1(t),0\le t\le
T_1(H_1))$ and $(\xi_2(t),0\le t\le T_2(H_2))$ is $\nu$ instead of
$\mu$. Fix an $(\F^2_t)$-stopping time $\bar t_2$ with $\bar t_2\le
T_2(H_2)$. From (\ref{SDE-xi}), (\ref{t1}), and Girsanov theorem
(\cite{RY}), there is an $(\F^1_t\times\F^2_{\bar t_2})$-Brownian
motion $\til B_1(t)$  such that $\xi_1(t_1)$ satisfies the
$(\F^1_{t_1}\times\F^2_{\bar t_2})$-adapted SDE for $0\le t_1\le
T_1(H_1)$: \BGE d\xi_1(t_1)=\sqrt\kappa d \til
B_1(t_1)+(3-\frac\kappa 2) \Big(\frac{A_{1,2}(t_1,\bar t_2)}{
A_{1,1}(t_1,\bar t_2)}+\frac{2A_{1,1}(t_1,\bar t_2)}{
A_{2,0}(t_1,\bar t_2)- A_{1,0}(t_1,\bar t_2)} \Big)d
t_1.\label{nu}\EDE From (\ref{pa-A_1-0}) and (\ref{nu}), we have
\BGE d A_{1,0}(t_1,\bar t_2)=A_{1,1}(t_1,\bar t_2)\sqrt\kappa d\til
B_1(t)+ \frac{(6-\kappa)A_{1,1}(t_1,\bar
t_2)^2dt_1}{A_{2,0}(t_1,\bar t_2)-A_{1,0}(t_1,\bar
t_2)}.\label{nu-pre-eta}\EDE Recall that $A_{1,0}(t_1,\bar
t_2)=\vphi_{2,t_1}(\bar t_2,\xi_1(t_1))=\xi_{1,\bar
t_2}(t_1)=\eta_{1,\bar t_2}(v_{1,\bar t_2}(t_1))$, and $ v_{1,\bar
t_2}'(t_1)=A_{1,1}(t_1,\bar t_2)^2$ (see (\ref{v'})). From
(\ref{nu-pre-eta}), there is a Brownian motion $\ha B_1(t_1)$ such
that \BGE d\eta_{1,\bar t_2}(s_1)=\sqrt\kappa d\ha
B_1(s_1)+\frac{(\kappa -6)ds_1} {\eta_{1,\bar t_2}(s_1)-
A_{2,0}(v_{1,\bar t_2}^{-1}(s_1),\bar t_2)}.\label{eta-A-1}\EDE
Since $A_{2,0}(v_{1,\bar t_2}^{-1}(s_1),\bar t_2)=\vphi_{1,\bar
t_2}(v_{1,\bar t_2}^{-1}(s_1),\xi_2(\bar t_2))=\psi_{1,\bar
t_2}(s_1,\xi_2(\bar t_2))$ and $\psi_{1,\bar t_2}(s,\cdot)$, $0\le
s\le v_{1,\bar t_2}(T_1(H_1))$, are chordal Loewner maps driven by
$\eta_{1,\bar t_2}(s)$, so the chordal Loewner hulls $L_{1,\bar
t_2}(s)$, $0\le s\le v_{1,\bar t_2}(T_1(H_1))$, driven by
$\eta_{1,\bar t_2}(s)$ is a part of chordal SLE$(\kappa, \kappa-6)$
process started from $\eta_{1,\bar t_2}(0)=\vphi_2(\bar t_2,x_1)$
with force point at $A_{2,0}(v_{1,\bar t_2}^{-1}(0),\bar
t_2)=\xi_2(\bar t_2)$. Thus after a time-change, it is a chordal
SLE$_\kappa$ in $\HH$ from $\vphi_2(\bar t_2,x_1)$ to $\xi_2(\bar
t_2)$. Note that $\vphi_2(\bar t_2,\cdot)^{-1}$ maps $\HH$
conformally onto $\HH\sem\beta_2((0,\bar t_2])$, maps $L_{1,\bar
t_2}(v_{1,\bar t_2}(t_1))$ onto $K_1(t_1)=\beta_1((0,t_1])$, and
takes $\vphi_2(\bar t_2,x_1)$ and $\xi_2(\bar t_2)$ to $x_1$ and
$\beta_2(\bar t_2)$, respectively. Thus $\beta_1(t)$, $0\le t\le
T_1(H_1)$, is the time-change of a chordal SLE$_\kappa$ trace in
$\HH\sem\beta_2((0,\bar t_2])$ from $x_1$ to $\beta_2(\bar t_2)$,
stopped on hitting $\lin{\HH\sem H_1}$. Similarly, for any
$(\F^1_t)$-stopping time $\bar t_1$ with $\bar t_1\le T_1(H_1)$,
$\beta_2(t)$, $0\le t\le T_2(H_2)$, is a time-change of a chordal
SLE$_\kappa$ trace in $\HH\sem \beta_1((0,\bar t_1])$ from $x_2$ to
$\beta_1(\bar t_1)$ stopped on hitting $\lin{\HH\sem H_2}$.

\section{Constructing New Martingales}\label{Const}
\begin{Theorem} For any $(H_1^m,H_2^m)\in\HP$, $1\le m\le n$, there
is a continuous function $M_*(t_1,t_2)$ defined on $[0,\infty]^2$
that satisfies the following properties: (i) $M_*=M$ on
$[0,T_1(H_1^m)]\times[0,T_2(H_2^m)]$ for $m=1,\dots,n$; (ii)
$M_*(t,0)=M_*(0,t)=1$ for any $t\ge 0$; (iii) $M_*(t_1,t_2)\in
[C_1,C_2]$ for any $t_1,t_2\ge 0$, where $C_2>C_1>0$ are constants
depending only on $H_j^m$, $j=1,2$, $1\le m\le n$; (iv) for any
$(\F^2_t)$-stopping time $\bar t_2$, $(M_*(t_1,\bar t_2),t_1\ge 0)$
is a bounded continuous $(\F^1_{t_1}\times\F^2_{\bar t_2})_{t_1\ge
0}$-martingale; and (v) for any $(\F^1_t)$-stopping time $\bar t_1$,
$(M_*(\bar t_1, t_2),t_2\ge 0)$ is a bounded continuous $(\F^1_{\bar
t_1}\times\F^2_{t_2})_{t_2\ge 0}$-martingale.
 \label{martg}
\end{Theorem}
{\bf Proof.} We will first define $M_*$ and then check its
properties. The first quadrant $[0,\infty]^2$ is divided by the
horizontal or vertical lines $\{x_j=T_j(H_j^m)\}$, $1\le m\le n$,
$j=1,2$, into small rectangles, and $M_*$ is piecewise defined on
each rectangle. Theorem \ref{M-I} will be used to prove the
martingale properties.

Let $\N_n:=\{k\in\N:k\le n\}$. Write $T_j^k$ for $T_j(H^k_j)$,
$k\in\N_n$, $j=1,2$. Let $S\subset\N_n$ be such that $\cup_{k\in
S}[0,T_1^k]\times[0,T_2^k]=\cup_{k=1}^n [0,T_1^k]\times[0,T_2^k]$,
and $\sum_{k\in S} k\le \sum_{k\in S'}k$ if $S'\subset\N_n$ also
satisfies this property. Such $S$ is a random nonempty set, and
$|S|\in\N_n$ is a random number. Define an partial order
``$\preceq$'' on $[0,\infty]^2$ such that
$(s_1,s_2)\preceq(t_1,t_2)$ iff $s_1\le t_1$ and $s_2\le t_2$. If
$(s_1,s_2)\preceq(t_1,t_2)$ and $(s_t,s_2)\ne (t_1,t_2)$,  we write
$(s_1,s_2)\prec(t_1,t_2)$. Then for each $m\in\N_n$ there is $k\in
S$ such that $(T_1^m,T_2^m)\preceq(T_1^k,T_2^k)$; and for each $k\in
S$ there is no $m\in\N_n$ such that
$(T_1^k,T_2^k)\prec(T_1^m,T_2^m)$.

There is a map $\sigma$ from $\{1,\dots,|S|\}$ onto $S$ such that if
$1\le k_1<k_2\le |S|$, then \BGE
T_1^{\sigma(k_1)}<T_1^{\sigma(k_2)}, \quad
T_2^{\sigma(k_1)}>T_2^{\sigma(k_2)}.\label{<}\EDE Define
$T_1^{\sigma(0)}=T_2^{\sigma(|S|+1)}=0$ and
$T_1^{\sigma(|S|+1)}=T_2^{\sigma(0)}=\infty$. Then (\ref{<}) still
holds for $0\le k_1<k_2\le |S|+1$.

Extend the definition of $M$ to
$[0,\infty]\times\{0\}\cup\{0\}\times[0,\infty]$ such that
$M(t,0)=M(0,t)=1$ for $t\ge 0$. Fix $(t_1,t_2)\in[0,\infty]^2$.
There are $k_1\in \N_{|S|+1}$ and $k_2\in\N_{|S|}\cup\{0\}$ such
that \BGE T^{\sigma(k_1-1)}_1\le t_1\le T^{\sigma(k_1)}_1,\quad
T^{\sigma(k_2+1)}_2\le t_2\le T^{\sigma(k_2)}_2.\label{k1k2}\EDE If
$k_1\le k_2$, let \BGE M_*(t_1,t_2)=M(t_1,t_2).\label{M*M}\EDE It
$k_1\ge k_2+1$, let \BGE
M_*(t_1,t_2)=\frac{M(T_1^{\sigma(k_2)},t_2)M(T_1^{\sigma(k_2+1)},
T_2^{\sigma(k_2+1)})\cdots
M(T_1^{\sigma(k_1-1)},T_2^{\sigma(k_1-1)})M(t_1,T_2^{\sigma(k_1)}) }
{M(T_1^{\sigma(k_2)},T_2^{\sigma(k_2+1)})\cdots
M(T_1^{\sigma(k_1-2)},T_2^{\sigma(k_1-1)})
M(T_1^{\sigma(k_1-1)},T_2^{\sigma(k_1)}) }\label{M*}\EDE In the
above formula, there are $k_1-k_2+1$ terms in the numerator, and
$k_1-k_2$ terms in the denominator. For example, if $k_1-k_2=1$,
then
$$M_*(t_1,t_2)=M(T_1^{\sigma(k_2)},t_2)M(t_1,T_2^{\sigma(k_1)})
/M(T_1^{\sigma(k_2)},T_2^{\sigma(k_1)}).$$

We need to show that $M_*(t_1,t_2)$ is well-defined. First, we show
that the $M(\cdot,\cdot)$ in (\ref{M*M}) and (\ref{M*}) are defined.
Note that $M$ is defined on
$$Z:=\bigcup_{k=0}^{|S|+1}[0,T_1^{\sigma(k)}]\times[0,T_2^{\sigma(k)}].$$
If $k_1\le k_2$ then $t_1\le T_1^{\sigma(k_1)}\le T_1^{\sigma(k_2)}$
and $t_2\le T_2^{\sigma(k_2)}$, so $(t_1,t_2)\in Z$. Thus
$M(t_1,t_2)$ in (\ref{M*M}) is defined. Now suppose $k_1\ge k_2+1$.
Since $t_2\le T_2^{\sigma(k_2)}$ and $t_1\le T_1^{\sigma(k_1)}$, so
$(T_1^{\sigma(k_2)},t_2),(t_1,T_2^{\sigma(k_1)})\in Z$. It is clear
that $(T_1^{\sigma(k)}, T_2^{\sigma(k)})\in Z$ for $k_2+1\le k\le
k_1-1$. Thus the $M(\cdot,\cdot)$ in the numerator of (\ref{M*}) are
defined. For $k_2\le k\le k_1-1$, $T_1^{\sigma(k)}\le
T_1^{\sigma(k+1)}$, so $(T_1^{\sigma(k)}, T_2^{\sigma(k+1)})\in Z$.
Thus the $M(\cdot,\cdot)$ in the denominator of (\ref{M*}) are
defined.

Second, we show that the value of $M_*(t_1,t_2)$ does not depend on
the choice of $(k_1,k_2)$ that satisfies (\ref{k1k2}). Suppose
(\ref{k1k2}) holds with $(k_1,k_2)$ replaced by $(k_1',k_2)$, and
$k_1'\ne k_1$. Then $|k_1'-k_1|=1$. We may assume $k_1'=k_1+1$. Then
$t_1=T_1^{\sigma(k_1)}$. Let $M_*'(t_1,t_2)$ denote the
$M_*(t_1,t_2)$ defined using $(k_1',k_2)$. There are three cases.
Case 1. $k_1<k_1'\le k_2$. Then from (\ref{M*M}),
$M_*'(t_1,t_2)=M(t_1,t_2)=M_*(t_1,t_2)$. Case 2. $k_1=k_2$ and
$k_1'-k_2=1$. Then $T_1^{\sigma(k_2)}=T_1^{\sigma(k_1)}=t_1$. So
from (\ref{M*M}) and (\ref{M*}),
$$M_*'(t_1,t_2)={M(T_1^{\sigma(k_2)},t_2)M(t_1,T_2^{\sigma(k_1)})}/
{M(T_1^{\sigma(k_2)},T_2^{\sigma(k_1)})}=M(t_1,t_2)=M_*(t_1,t_2).$$
Case 3. $k_1'>k_1>k_2$. From (\ref{M*}) and that
$T_1^{\sigma(k_1)}=t_1$, we have
$$M_*'(t_1,t_2)=\frac{M(T_1^{\sigma(k_2)},t_2)M(T_1^{\sigma(k_2+1)},T_2^{\sigma(k_2+1)})
\cdots
M(T_1^{\sigma(k_1)},T_2^{\sigma(k_1)})M(t_1,T_2^{\sigma(k_1+1)})}
{M(T_1^{\sigma(k_2)},T_2^{\sigma(k_2+1)})\cdots
M(T_1^{\sigma(k_1-1)},T_2^{\sigma(k_1)})
M(T_1^{\sigma(k_1)},T_2^{\sigma(k_1+1)})}$$
$$=\frac{M(T_1^{\sigma(k_2)},t_2)M(T_1^{\sigma(k_2+1)},T_2^{\sigma(k_2+1)})
\cdots M(t_1,T_2^{\sigma(k_1)})}
{M(T_1^{\sigma(k_2)},T_2^{\sigma(k_2+1)})\cdots
M(T_1^{\sigma(k_1-1)},T_2^{\sigma(k_1)})}=M_*(t_1,t_2).$$ Similarly,
if (\ref{k1k2}) holds with $(k_1,k_2)$ replaced by $(k_1,k_2')$,
then $M_*(t_1,t_2)$ defined using $(k_1,k_2')$ has the same value as
$M(t_1,t_2)$. Thus $M_*$ is well-defined.

From the definition, it is clear that for each $k_1\in \N_{|S|+1}$
and $k_2\in\N_{|S|}\cup\{0\}$, $M_*$ is continuous on
$[T_1^{\sigma(k_1-1)},T_1^{\sigma(k_1)}]
\times[T_2^{\sigma(k_2+1)},T_1^{\sigma(k_2)}]$. Thus $M_*$ is
continuous on $[0,\infty]^2$. Let $(t_1,t_2)\in[0,\infty]^2$.
Suppose $(t_1,t_2)\in[0,T_1^m]\times[0,T_2^{m}]$ for some
$m\in\N_n$. There is $k\in \N_{|S|}$ such that
$(T_1^m,T_2^m)\preceq(T_1^{\sigma(k)},T_2^{\sigma(k)})$. Then we may
choose $k_1\le k$ and $k_2\ge k$ such that (\ref{k1k2}) holds, so
$M_*(t_1,t_2)=M(t_1,t_2)$. Thus (i) is satisfied. If $t_1=0$, we may
choose $k_1=1$ in (\ref{k1k2}). Then either $k_1\le k_2$ or $k_2=0$.
If $k_1\le k_2$ then $M_*(t_1,t_2)=M(t_1,t_2)=1$ because $t_1=0$. If
$k_2=0$, then
$$M_*(t_1,t_2)=M(T_1^{\sigma(0)},t_2)M(t_1,T_2^{\sigma(1)})
/M(T_1^{\sigma(0)},T_2^{\sigma(1)})=1$$ because
$T_1^{\sigma(0)}=t_1=0$. Similarly, $M_*(t_1,t_2)=0$ if $t_2=0$. So
(ii) is also satisfied. And (iii) follows from Lemma \ref{bound} and
the definition of $M_*$.

\vskip 3mm

Now we prove (iv). Suppose $(t_1,t_2)\in[0,\infty]^2$ and $t_2\ge
\vee_{m=1}^n T_2^m=T_2^{\sigma(1)}$. Then (\ref{k1k2}) holds with
$k_2=0$ and some $k_1\in\{1,\dots,|S|+1\}$. So $k_1\ge k_2+1$. Since
$T_1^{\sigma(k_2)}=0$ and $M(0,t)=1$ for any $t\ge 0$, so from
(\ref{M*}) we have
$$M_*(t_1,t_2)=\frac{M(T_1^{\sigma(k_2+1)},
T_2^{\sigma(k_2+1)})\cdots
M(T_1^{\sigma(k_1-1)},T_2^{\sigma(k_1-1)})M(t_1,T_2^{\sigma(k_1)}) }
{M(T_1^{\sigma(k_2+1)},T_2^{\sigma(k_2+2)})\cdots
M(T_1^{\sigma(k_1-1)},T_2^{\sigma(k_1)}) }.$$ The right-hand side of
the above equality has no $t_2$. So
$M_*(t_1,t_2)=M_*(t_1,\vee_{m=1}^n T_2^m)$ for any $t_2\ge
\vee_{m=1}^n T_2^m$. Similarly, $M_*(t_1,t_2)=M_*(\vee_{m=1}^n
T_1^m,t_2)$ for any $t_1\ge \vee_{m=1}^n T_1^m$.

 Fix an $(\F^2_t)$-stopping time $\bar t_2$. Since
$M_*(\cdot,\bar t_2)=M_*(\cdot,\bar t_2\wedge(\vee_{m=1}^n T_2^m))$,
and $\bar t_2\wedge(\vee_{m=1}^n T_2^m)$ is also an
$(\F^2_t)$-stopping time, so we may assume that $\bar t_2\le
\vee_{m=1}^n T_2^m$. Let $I_0=\bar t_2$. For $s\in \N\cup\{0\}$,
define \BGE R_s=\sup\{T_1^m:m\in\N_n,T_2^m\ge I_s\};\quad
I_{s+1}=\sup\{T_2^m:m\in\N_n,T_2^m<I_s,T_1^m>R_s\}.\label{RI}\EDE
Here we set $\sup(\emptyset)=0$. Then we have a non-decreasing
sequence $(R_s)$ and a non-increasing sequence $(I_s)$. Let $S$ and
$\sigma(k)$, $0\le k\le |S|+1$, be as in the definition of $M_*$.
From the property of $S$, for any $s\in\N\cup\{0\}$, \BGE
R_s=\sup\{T_1^k:k\in S,T_2^k\ge I_s\}.\label{R}\EDE Suppose for some
$s\in\N\cup\{0\}$, there is $m\in\N_n$ that satisfies $T_2^m<I_s$
and $T_1^m>R_s$. Then there is $k\in S$ such that $T_j^k\ge T_j^m$,
$j=1,2$. If $T_2^k\ge I_s$, then from (\ref{R}) we have $R_s\ge
T_1^k\ge T_1^m$, which contradicts that $T_1^m>R_s$. Thus
$T_2^k<I_s$. Now $T_2^k<I_s$, $T_1^k\ge T_1^m>R_s$, and $T_2^k\ge
T_2^m$. Thus for any $s\in\N\cup\{0\}$, \BGE
I_{s+1}=\sup\{T_2^k:k\in S,T_2^k<I_s,T_1^k>R_s\}.\label{I}\EDE

First suppose $\bar t_2>0$. Since $\bar t_2\le \vee_{m=1}^n
T_2^m=T_2^{\sigma(0)}$, so there is a unique $k_2\in \N_{|S|}$ such
that $T_2^{\sigma(k_2)}\ge \bar t_2>T_2^{\sigma(k_2+1)}$. From
(\ref{R}) and (\ref{I}), we have $R_s=T_1^{\sigma(k_2+s)}$ for $0\le
s\le |S|-k_2$; $R_s=T_1^{\sigma(|S|)}$ for $s\ge |S|-k_2$;
$I_s=T_2^{\sigma(k_2+s)}$ for $1\le s\le |S|-k_2$; and $I_s=0$ for
$s\ge |S|-k_2+1$. Since $R_0= T_1^{\sigma(k_2)}$ and $\bar t_2\le
T_2^{\sigma(k_2)}$, so from (i),
 \BGE M_*(t_1,\bar
t_2)=M(t_1,\bar t_2),\quad\mbox{for } t_1\in[0,R_0].\label{0,R0}\EDE
Suppose $t_1\in [R_{s-1},R_{s}]$ for some $s\in\N_{|S|-k_2}$. Let
$k_1=k_2+s$. Then $T_1^{\sigma(k_1-1)}\le t_1\le T_1^{\sigma(k_1)}$.
Since $I_s=T_2^{\sigma(k_2+s)}=T_2^{\sigma(k_1)}$, so from
(\ref{M*}), \BGE M_*(t_1,\bar t_2)/M_*(R_{s-1},\bar
t_2)=M(t_1,I_{s})/M(R_{s-1}, I_{s}),\quad\mbox{for }t_1\in
[R_{s-1},R_{s}].\label{M/M-1}\EDE Note that if $s\ge |S|-k_2+1$,
(\ref{M/M-1}) still holds because $R_s=R_{s-1}$. Suppose $t_1\ge
R_{n}$. Since $n\ge |S|-k_2$, so $R_n=T_1^{\sigma(|S|)}=\vee_{m=1}^n
T_1^m$. From the discussion at the beginning of the proof of (iv),
we have \BGE M_*(t_1,\bar t_2)=M_*(R_n,\bar t_2),\quad\mbox{for
}t_1\in[R_n,\infty].\label{M=M}\EDE If $\bar t_2=0$,
(\ref{0,R0}$\sim$\ref{M=M}) still hold because all $I_s=0$ and so
$M_*(t_1,\bar t_2)=M(t_1,I_s)=M(t_1,0)=1$ for any $t_1\ge 0$.

Let $R_{-1}=0$. We claim that for each $s\in\N\cup\{0\}$, $R_s$ is
an $(\F^1_t\times\F^2_{\bar t_2})_{t\ge 0}$-stopping time, and
$I_{s}$ is $\F_{R_{s-1},\bar t_2}$-measurable. Recall that
$\F_{R_{s-1},\bar t_2}$ is the $\sigma$-algebra obtained from the
filtration $(\F^1_t\times\F^2_{\bar t_2})_{t\ge 0}$ and its stopping
time $R_{s-1}$. It is clear that $R_{-1}=0$ is an
$(\F^1_t\times\F^2_{\bar t_2})_{t\ge 0}$-stopping time, and
$I_0=\bar t_2$ is $\F_{R_{-1},\bar t_2}$-measurable. Now suppose
$I_{s}$ is $\F_{R_{s-1},\bar t_2}$-measurable. Since $I_s\le \bar
t_2$ and $R_{s-1}\le R_s$, so for any $t\ge 0$, $\{R_s\le
t\}=\{R_{s-1}\le t\}\cap {\cal E}_t$, where
$${\cal E}_t=\bigcap_{m=1}^n(\{T_2^m<I_s\}\cup\{T_1^m\le t\})=
 \bigcap_{m=1}^n(\cup_{q\in\Q}(\{T_2^m<q\le \bar t_2\}\cap\{q<I_s\})\cup \{T_1^m\le t\} ).$$
Thus ${\cal E}_t\in \F_{R_{s-1},\bar t_2}\vee
(\F^1_t\times\F^2_{\bar t_2})$, and so $\{R_s\le t\}\in
\F^1_t\times\F^2_{\bar t_2}$ for any $t\ge 0$. Therefore $R_s$ is an
$(\F^1_t\times\F^2_{\bar t_2})_{t\ge 0}$-stopping time. Next we
consider $I_{s+1}$. For any $h\ge 0$, $$\{I_{s+1}>
h\}=\cup_{m=1}^n(\{h<T_2^m<I_s\}\cap\{T_1^m>R_s\})$$
$$=\cup_{m=1}^n(\cup_{q\in\Q}(\{h<T_2^m<q< \bar t_2\}\cap\{q<I_s\})
\cap\{T_1^m>R_s\})\in \F_{R_{s},\bar t_2}.$$ Thus $I_{s+1}$ is
$\F_{R_{s},\bar t_2}$-measurable. So the claim is proved by
induction.

Since $\bar t_2\le\vee_{m=1}^n T_2^m<T_2$, so from Theorem
\ref{M-I}, for any $s\in\N_n$, $(M(R_{s-1}+t,I_s),0\le
t<T_1(I_s)-R_{s-1})$ is a continuous $(\F_{R_{s-1}+t,\bar
t_2})_{t\ge 0}$-local martingale. For $m\in\N_n$, if $T_2^m\ge I_s$,
then $T_1^m<T_1(T_2^m)\le T_1(I_s)$. So from (\ref{RI}) we have
$R_s<T_1(I_s)$. From (\ref{M/M-1}), we find that
 $(M_*(R_{s-1}+t,\bar t_2),0\le t\le R_s-R_{s-1})$ is
a continuous $(\F_{R_{s-1}+t,\bar t_2})_{t\ge 0}$-local martingale
for any $s\in\N_n$. From Theorem \ref{M} and (\ref{0,R0}),
$(M_*(t,\bar t_2),0\le t\le R_0)$ is a continuous $(\F_{t,\bar
t_2})_{t\ge 0}$-local martingale. From (\ref{M=M}), $(M_*(R_n+t,\bar
t_2),t\ge 0)$ is a continuous $(\F_{R_n+t,\bar t_2})_{t\ge 0}$-local
martingale. Thus $(M_*(t,\bar t_2),t\ge 0)$ is a continuous
$(\F_{t,\bar t_2})_{t\ge 0}$-local martingale. Since  by (iii)
$M_*(t_1,t_2)\in[C_1, C_2]$, so this local martingale is a bounded
martingale. Thus (iv) is satisfied. Finally, (v) follows from the
symmetry in the definition (\ref{M*M})
 and (\ref{M*}) of $M_*$. $\Box$

\section{Coupling Measures} \label{coup}
{\bf Proof of Theorem \ref{Theorem}.} From conformal invariance, we
may assume that $D=\HH$, $a=x_1$ and $b=x_2$. Let $\xi_j(t)$ and
$\beta_j(t)$, $0\le t\le T_j$, $j=1,2$, be as in Section
\ref{Two-Dim}. For $j=1,2$, let $\mu_j$ denote the distribution of
$(\xi_j(t),0\le t\le T_j)$. Let $\mu=\mu_1\times \mu_2$. Then $\mu$
is the joint distribution of $\xi_1$ and $\xi_2$, since they are
independent.

Let $\ha\C=\C\cup\{\infty\}$ be the Riemann sphere with spherical
metric. Let $\Gamma_{\ha\C}$ denote the space of nonempty compact
subsets of $\ha\C$ endowed with Hausdorff metric. Then
$\Gamma_{\ha\C}$ is a compact metric space. For a chordal Loewner
trace $\beta(t)$, $0\le t\le T$, let $G(\beta):=\{\beta(t):0\le t\le
T\}\in \Gamma_{\ha \C}$. For $j=1,2$, let $\bar\mu_j$ denote the
distribution of $G(\beta_j)$, which is a probability measure on
$\Gamma_{\ha\C}$. We want to prove that $\bar\mu_1=\bar\mu_2$. Let
$\bar\mu=\bar\mu_1\times\bar\mu_2$, which is the joint distribution
of $G(\beta_1)$ and $G(\beta_2)$.

Let $\HP_*$ be the set of $(H_1,H_2)\in\HP$ such that for $j=1,2$,
$H_j$ is a polygon whose vertices have rational coordinates. Then
$\HP_*$ is countable. Let $(H_1^m,H_2^m)$, $m\in\N$, be an
enumeration of $\HP_*$. For each $n\in\N$, let $M_*^n(t_1,t_2)$ be
the $M_*(t_1,t_2)$ given by Theorem \ref{martg} for $(H_1^m,H_2^m)$,
$1\le m\le n$, in the above enumeration.

For each $n\in\N$ define $\nu^n=(\nu^n_1,\nu^n_2)$ such that
${d\nu^n}/{d\mu}=M_*^n(\infty,\infty)$. From Theorem \ref{martg},
$M_*^n(\infty,\infty)>0$ and $\int
M_*^n(\infty,\infty)d\mu=\EE[M_*^n(\infty,\infty)]=1$, so $\nu^n$ is
a probability measure. Then
$d\nu^n_1/d\mu_1=\EE[M_*^n(\infty,\infty)|\F^2_\infty]=M_*^n(\infty,0)=1$.
Thus $\nu^n_1=\mu_1$. Similarly, $\nu^n_2=\mu_2$. So each $\nu^n$ is
a coupling of $\mu_1$ and $\mu_2$.

For each $n\in\N$, suppose $(\zeta_1^n(t),0\le t\le S^n_1)$ and
$(\zeta_2^n(t),0\le t\le S^n_2)$ have the joint distribution
$\nu^n$. Let $\gamma^n_j(t)$, $0\le t\le S_j$, $j=1,2$, be the
chordal Loewner trace driven by $\zeta^n_j$. Let $\bar
\nu^n=(\bar\nu^n_1,\bar\nu^n_2)$ denote the joint distribution of
$G(\gamma^n_1)$ and $G(\gamma^n_2)$. Since $\Gamma_{\ha\C}\times
\Gamma_{\ha\C}$ is compact, so $(\bar\nu^n,n\in\N)$ has a
subsequence $(\bar\nu^{n_k}:k\in\N)$ that converges weakly to some
probability measure $\bar\nu=(\bar\nu_1,\bar\nu_2)$ on
$\Gamma_{\ha\C}\times \Gamma_{\ha\C}$. Then for $j=1,2$,
$\bar\nu^{n_k}_j\to\bar\nu_j$ weakly. For $n\in\N$ and $j=1,2$,
since $\nu^n_j=\mu_j$, so $\bar \nu^n_j=\bar \mu_j$. Thus
$\bar\nu_j=\bar\mu_j$, $j=1,2$. So $\bar\nu_j$, $j=1,2$, is
supported by the space of graphs of crosscuts in $\HH$. From
Proposition \ref{chordal-Loewner-chain}, there are $\zeta_1\in
C([0,S_1])$ and $\zeta_2\in C([0,S_2])$ such that the joint
distribution of $G(\gamma_1)$ and $G(\gamma_2)$ is $\bar\nu$, where
$\gamma_j(t)$ is the chordal Loewner trace driven by $\zeta_j(t)$,
$j=1,2$.

Now fix $m\in\N$. From Theorem \ref{M}, $M(T_1(H_1^m),T_2(H_2^m))$
is positive and $\F^1_{T_1(H_1^m)}\times
\F^2_{T_2(H_2^m)}$-measurable, and $\int
M(T_1(H_1^m),T_2(H_2^m))d\mu=1$. Define $\nu_{(m)}$ on
$\F^1_{T_1(H_1^m)}\times \F^2_{T_2(H_2^m)}$ such that
$d\nu_{(m)}/d\mu=M(T_1(H_1^m),T_2(H_2^m))$. Then $\nu_{(m)}$ is a
probability measure. From Theorem \ref{martg}, if $n\ge m$, then
$$\frac{d\nu^n}{d\mu}\Big|_{\F^1_{T_1(H_1^m)}\times \F^2_{T_2(H_2^m)}}
=\EE[M_*^n(\infty,\infty)|\F^1_{T_1(H_1^m)}\times \F^2_{T_2(H_2^m)}]
$$$$=M_*^n(T_1(H_1^m),T_2(H_2^m))=M(T_1(H_1^m),T_2(H_2^m)).$$ Thus
$\nu_{(m)}$ equals to the restriction of $\nu^n$ to
$\F^1_{T_1(H_1^m)}\times \F^2_{T_2(H_2^m)}$ if $n\ge m$.

For a chordal Loewner trace $\gamma(t)$, $0\le t\le S$, and a hull
$H$ in $\HH$ w.r.t.\ $0$ that contains some neighborhood of
$\gamma(0)$ in $\HH$, let $G_H(\gamma):=\{\gamma(t):0\le t\le
T_H\}\in \Gamma_{\ha\C}$, where $T_H$ is the first $t$ such that
$\gamma(t)\in\lin{\HH\sem H}$ or $t=S$. Then $G_H(\gamma)\subset
G(\gamma)$. Let $\bar\nu^n_{(m)}$ denote the distribution of
$(G_{H_1^m}(\gamma^{n}_1),G_{H_2^m}(\gamma^{n}_2))$. Then $\bar
\nu^n_{(m)}$ is determined by the distribution of
$(\zeta^{n}_1,\zeta^{n}_2)$ restricted to $\F^1_{T_1(H_1^m)}\times
\F^2_{T_2(H_2^m)}$, which equals to $\nu_{(m)}$ if $n\ge m$. Let
$\bar\nu_{(m)}=\bar\nu^m_{(m)}$. Then
$\bar\nu^n_{(m)}=\bar\nu_{(m)}$ for $n\ge m$.

Let $\tau^{n_k}_{(m)}$ denote the distribution of
$(G(\gamma^{n_k}_1),G(\gamma^{n_k}_2),G_{H_1^m}(\gamma^{n_k}_1),G_{H_2^m}(\gamma^{n_k}_2)
)$. Then $\tau^{n_k}_{(m)}$ is supported by $\Xi$, which is the set
of $(L_1,L_2,F_1,F_2)\in \Gamma_{\ha\C}^4$ such that $F_j\subset
L_j$ for $j=1,2$. It is easy to check that $\Xi$ is a closed subset
of $\Gamma_{\ha\C}^4$. Then $(n_k)$ has a subsequence $(n'_k)$ such
that $\tau^{n'_k}_{(m)}$ converges weakly to some probability
measure $\tau_{(m)}$ on $\Xi$. Since the marginal of
$\tau^{n'_k}_{(m)}$ at the first two variables equals to $\bar
\nu^{n'_k}$, and $\bar \nu^{n'_k}\to\bar\nu$ weakly, so the marginal
of $\tau_{(m)}$ at the first two variables equals to $\bar\nu$.
Since the marginal of $\tau^{n'_k}_{(m)}$ at the last two variables
equals to $\bar \nu^{n'_k}_{(m)}$, which equals to $\bar\nu_{(m)}$
if $n'_k\ge m$, so the marginal of $\tau_{(m)}$ at the last two
variables equals to $\bar\nu_{(m)}$.

Let the $\Xi$-valued random variable $(L_1,L_2,F_1,F_2)$ has the
distribution $\tau_{(m)}$. Then $\bar\nu$ is the distribution of
$(L_1,L_2)$ and $\bar\nu_{(m)}$ is the distribution of $(F_1,F_2)$.
Note that $\bar\nu_{(m)}$ is supported by the space of pairs of
curves $(\alpha_1,\alpha_2)$ such that for $j=1,2$, $\alpha_j$ is a
simple curve whose one end is $x_j$, the other end lies on $\pa
H_j^m\cap\HH$, and whose other part lies in the interior of $H_j^m$.
For $j=1,2$, since $L_j=G(\gamma_j)$, so from the properties of
$\Xi$ and $\bar\nu_{(m)}$, we have $F_j=G_{H_j^m}(\gamma_j)$, which
means that $(G_{H_1^m}(\gamma_1),G_{H_2^m}(\gamma_2))$ has the
distribution $\bar\nu_{(m)}$. Since the distribution of
$(G_{H_1^m}(\gamma_1),G_{H_2^m}(\gamma_2))$ determines the
distribution of $(\zeta_1,\zeta_2)$ restricted to
$\F^1_{T_1(H_1^m)}\times \F^2_{T_2(H_2^m)}$, so the the distribution
of $(\zeta_1,\zeta_2)$ restricted to $\F^1_{T_1(H_1^m)}\times
\F^2_{T_2(H_2^m)}$ equals to $\nu_{(m)}$. Since
$d\nu_{(m)}/d\mu=M(T_1(H_1^m),T_2(H_2^m))$, so from the discussion
after the proof of Theorem \ref{bound}, for any $(\F^2_t)$-stopping
time $\bar t_2$ with $\bar t_2\le T_2(H_2^m)$, $(\gamma_1(t),0\le
t\le T_1(H^m_1))$ is a time-change of a chordal SLE$_\kappa$ trace
in $\HH\sem \gamma_2((0,\bar t_2])$ from $x_1$ to $\gamma_2(\bar
t_2)$ stopped on hitting $\lin{\HH\sem H^m_1}$.

Now fix an $(\F^2_t)$-stopping time $\bar t_2$ with $\bar t_2<T_2$.
Recall that $T_1(\bar t_2)$ is the maximal such that
$\gamma_1([0,T_1(\bar t_2)))$ is disjoint from $\gamma_2([0,\bar
t_2])$. For $n\in\N$, define
$$R_n=\sup\{T_1(H^m_1):m\in\N_n,\bar t_2\le T_2(H^m_2)\}.$$
Here we set $\sup(\emptyset)=0$. Then for any $t\ge 0$,
$$\{R_n\le t\}=\cap_{m=1}^n(\{\bar t_2>
T_2(H^m_2)\}\cup\{T_1(H^m_1)\le t\})\in\F^1_t\times\F^2_{\bar
t_2}.$$ So $R_n$ is an $(\F^1_t\times\F^2_{\bar t_2})_{t\ge
0}$-stopping time for each $n\in\N$. For $m\in\N_n$, let $\bar
t_2^m=\bar t_2\wedge T_2(H^m_2)$.  Then $\bar t_2^m$ is an
$(\F^2_t)$-stopping time, and $\bar t_2^m\le T_2(H^m_2)$. From the
last paragraph, we conclude that $\gamma_1(t)$, $0\le t\le
T_1(H^m_1)$, is a time-change of a part of chordal SLE$_\kappa$
trace in $\HH\sem \gamma_1((0,\bar t_2^m])$ from $x_1$ to
$\gamma_2(\bar t_2^m)$. Let ${\cal E}_{n,m}=\{\bar t_2\le
T_2(H^m_2)\}\cap\{R_n=T_1(H^m_1)\}$. Since on each ${\cal E}_{n,m}$,
$\bar t_2=\bar t_2^m$ and $R_n=T_1(H^m_1)$, and
$\{R_n>0\}=\cup_{m=1}^n{\cal E}_{n,m}$, so $\gamma_1(t)$, $0\le t\le
R_n$, is a time-change of a part of chordal SLE$_\kappa$ trace in
$\HH\sem \gamma_1((0,\bar t_2])$ from $x_1$ to $\gamma_2(\bar t_2)$.
Let $R_\infty=\vee_{n=1}^\infty R_n$. Then $\gamma_1(t)$, $0\le t<
R_\infty$, is a time-change of a part of chordal SLE$_\kappa$ trace
in $\HH\sem \gamma_1((0,\bar t_2])$ from $x_1$ to $\gamma_2(\bar
t_2)$.

For each $n\in\N$ and $m\in\N_n$, if $\bar t_2\le T_2(H^m_2)$ then
$T_1(H^m_2)<T_1(\bar t_2)$, so $R_n<T_1(\bar t_2)$. Thus
$R_\infty\le T_1(\bar t_2)$. If $R_\infty<T_1(\bar t_2)$, then
$\gamma_1((0,R_\infty])$ is disjoint from $\gamma_2((0,\bar t_2])$,
so there is $(H^m_1,H^m_2)\in\HP_*$ such that
$\gamma_1((0,R_\infty])$ and $\gamma_2((0,\bar t_2])$ are contained
in the interiors of $H^m_1$ and $H^m_2$, respectively. Then $\bar
t_2\le T_2(H^m_2)$ and $R_m\le R_\infty<T_1(H^m_1)$, which
contradicts the definition of $R_m$. Thus $R_\infty=T_1(\bar t_2)$.
So $\gamma_1(t)$, $0\le t< T_1(\bar t_2)$, is a time-change of a
part of chordal SLE$_\kappa$ trace in $\HH\sem \gamma_1((0,\bar
t_2])$ from $x_1$ to $\gamma_2(\bar t_2)$. From the definition of
$T_1(\bar t_2)$ we have $\gamma_1(T_1(\bar t_2))\in G(\gamma_2)$.
Thus $\gamma_1(t)$, $0\le t< T_1(\bar t_2)$, is a time-change of a
full chordal SLE$_\kappa$ trace in $\HH\sem \gamma_1((0,\bar t_2])$
from $x_1$ to $\gamma_2(\bar t_2)$. Since $\kappa\in(0,4]$, so
almost surely $\gamma_1(T_1(\bar t_2))=\gamma_2(\bar t_2)$. Thus
$\gamma_2(\bar t_2)\in G(\gamma_1)$ almost surely.

For $n\in\N$ and $q\in\Q_{\ge 0}$, let $\bar t_2^{n,q}=q\wedge
T_2(H^n_2)$. Then each $\bar t_2^{n,q}$ is an $(\F^2_t)$-stopping
time with $\bar t_2^{q,n}<T_2$. Since $\N\times\Q_{\ge 0}$ is
countable, so almost surely $\gamma_2(\bar t_2^{q,n})\in
G(\gamma_1)$ for every $n\in\N$ and $q\in\Q_{\ge 0}$. Since $\Q_{\ge
0}$ is dense in $\R_{\ge 0}$, $\gamma_2$ is continuous, and
$G(\gamma_1)$ is closed, so almost surely for every $n\in \N$,
$\gamma_2([0,T_2(H^n_2)])\subset G(\gamma_1)$. Since
$T_2=\vee_{n=1}^\infty T_2(H^n_2)$, so $G(\gamma_2)\subset
G(\gamma_1)$ almost surely. Similarly, $G(\gamma_1)\subset
G(\gamma_2)$ almost surely. Thus $G(\gamma_1)= G(\gamma_2)$ almost
surely. Since for $j=1,2$, the distribution of $G(\gamma_j)$ equals
to the distribution of $G(\beta_j)$, which is the SLE$_\kappa$ trace
in $\HH$ from $x_j$ to $x_{3-j}$, so we are done. $\Box$

\end{document}